\documentclass{article}

\usepackage{amssymb}
\usepackage{amsmath}
\usepackage{mathtools}

\usepackage{tikz}
\usepackage{pict2e}
\usepackage[caption=false]{subfig}

\usepackage{booktabs}
\usepackage{makecell}
\usepackage{enumitem}

\usepackage{algorithm}
\usepackage[noend]{algpseudocode}

\usepackage[normalem]{ulem}
\usepackage[dvipsnames]{xcolor}
\usepackage[hidelinks]{hyperref}

\newcommand{\email}[1]{%
  \href{mailto:#1}{\texttt{#1}}%
}

\newenvironment{keywords}
  {\par\medskip\noindent\textbf{Keywords.}\ \ignorespaces}
  {\par\medskip}

\title{Residual-Driven Lifting Identification for Nonlinear-Manifold
Reduced-Order Models of Parametrized Linear PDEs}

\author{%
  Francesco A.\,B.\ Silva%
    \thanks{Department of Nuclear Engineering, Texas A\&M University,
            423 Spence St, College Station, TX 77843-3141, USA
            (\email{francesco.silva@tamu.edu},
             \email{jean.ragusa@tamu.edu}).}%
  \and
  Jean C.\ Ragusa\footnotemark[1]%
  \and
  Theron Guo%
    \thanks{Department of Mechanical Engineering, Massachusetts
            Institute of Technology, 77 Massachusetts Avenue,
            Cambridge, MA 02139, USA
            (\email{theron.guo@mit.edu}).}%
  \and
  Rudy Geelen%
    \thanks{Department of Aerospace Engineering, Texas A\&M University,
            701 H.R. Bright Building, College Station, TX 77843-3133,
            USA (\email{rudy.geelen@tamu.edu}).}%
}

\begin{document}

\maketitle
 
%% ------------------------------------------------------------
%% Abstract
%% ------------------------------------------------------------
\begin{abstract}
%% Text of abstract
We introduce a residual-driven procedure for training nonlinear-manifold reduced-order models for parametrized linear partial differential equations that, given prescribed latent and lifting spaces, identifies the nonlinear lifting without high-fidelity solution snapshots. The approximation is represented by a low-dimensional latent coordinate together with a nonlinear lifting into a richer reduced space. Rather than fitting the lifting to snapshot data, we determine it by minimizing a computable residual-based upper bound for the state error. For affinely parametrized operators, the resulting training objective admits an efficient offline--online decomposition, and the lifting update reduces to a sequence of low-dimensional weighted least-squares problems. Numerical evaluations on an advection--diffusion problem and a plane-strain elasticity benchmark show that the proposed approach substantially improves accuracy over linear subspaces. The resulting nonlinear models achieve accuracy comparable to snapshot-driven nonlinear-manifold training while avoiding high-fidelity snapshots in the lifting-identification stage.
\end{abstract}

\vspace{0.2cm}

\newpage
\section{Introduction}
\label{sec:Introduction}

Parametrized partial differential equations arise in many computational settings, including uncertainty quantification, design optimization, and real-time simulation. In these applications, repeatedly solving a high-fidelity discretization can be prohibitively expensive. Projection-based reduced-order models address this challenge by approximating the high-dimensional solution in a low-dimensional trial space, thereby replacing the full-order problem with a smaller reduced system~\cite{benner2017,  grepl2007, hesthaven2016, quarteroni2015}. Commonly, these reduced spaces are constructed using proper orthogonal decomposition (POD) or greedy reduced-basis methods, both of which rely on snapshots of the solution~\cite{grepl2005, Volkwein2008, Sirovich1987}. The accuracy of these methods, however, is limited by the approximation properties of the selected linear space. When the solution manifold has slowly decaying Kolmogorov widths, a linear subspace may require a large dimension to achieve the desired accuracy~\cite{peherstorfer2016, pinkus2012}.

Nonlinear manifold reduced-order models seek to overcome this limitation by approximating the solution set on a nonlinear manifold rather than in a fixed linear subspace. Related approaches include interpolation of reduced models on matrix manifolds~\cite{amsallem2008, amsallem2011}, autoencoder-based nonlinear manifolds~\cite{Carlberg2020}, and polynomial manifold approximations that enrich a low-dimensional latent representation~\cite{Barnett2022, Geelen2024, JAIN201780}. This added geometric flexibility can substantially improve accuracy. Most existing constructions, however, identify the nonlinear manifold from high-fidelity solution snapshots, so that the offline stage requires full-order model evaluations over a representative set of parameter values, together with the storage and processing of the resulting solution data.

An important exception is the construction of Jain and coauthors~\cite{JAIN201780}, which is snapshot-free in the sense that modal derivatives are used to define the manifold curvature directly from the system operators. The present work pursues a complementary snapshot-free strategy for lifting identification: rather than relying on derivative-based expansions, the nonlinear lifting is identified by minimizing a computable residual-based upper bound on the state error.

Specifically, this work considers nonlinear manifold approximations of the form
\begin{equation*}
    \mathbf{c}_\varepsilon(\boldsymbol{\mu}) = \mathbb{V}\mathbf{c}_\nu(\boldsymbol{\mu}) + \mathbb{U}\mathbb{X} \mathbf{h}\bigl(\mathbf{c}_\nu(\boldsymbol{\mu})\bigr),
\end{equation*}
where $\mathbf{c}_\nu(\boldsymbol{\mu})$ is a low-dimensional latent coordinate, $\mathbb{V}$ and $\mathbb{U}$ define prescribed latent and lifting spaces, $\mathbf{h}$ is a nonlinear feature map, and $\mathbb{X}$ is the lifting matrix to be inferred during training. In snapshot-driven approaches, $\mathbb{X}$ is determined by minimizing a loss function involving the discrepancy with respect to full-order solutions. In contrast, the proposed method identifies $\mathbb{X}$ using the governing operators, the right-hand side, and parameter samples, but without requiring full-order solution evaluations at those samples.

Throughout the paper, the term \emph{snapshot-free} is used in a precise sense: the identification of the nonlinear lifting matrix $\mathbb{X}$ does not use high-fidelity solution snapshots. Suitable latent and lifting spaces are assumed to be available from analytical considerations, spectral constructions, or previously computed bases. When the prescribed spaces are themselves built from snapshots, as for the proper orthogonal decomposition spaces used in one of our experiments, those snapshots enter the construction of the approximation spaces but not the identification of the lifting, which does not access the snapshot set. For brevity, subsequent occurrences of ``snapshot-free'' refer to this precise notion unless stated otherwise. The central question is therefore how to identify the nonlinear lifting without solution snapshots once the approximation spaces have been selected.

The proposed approach is motivated by residual-based \emph{a posteriori} error estimation. For coercive parametrized linear problems, the state error can be bounded by the dual norm of the residual scaled by the inverse coercivity constant. This observation allows the reconstruction-error objective, which depends on the unknown full-order solution, to be replaced by a residual-based objective that can be evaluated directly from the discrete operators. For affinely parametrized problems, this objective admits an efficient offline--online decomposition: after precomputing parameter-independent reduced blocks, the lifting update becomes a low-dimensional weighted least-squares problem that can be solved without assembling or storing solution snapshots.

The main contribution of this paper is a residual-driven training framework for polynomial-manifold reduced-order models that is snapshot-free with respect to lifting identification. The proposed approach identifies the nonlinear lifting matrix without high-fidelity solution snapshots, admits an efficient affine offline--online implementation, and is assessed on scalar advection--diffusion and plane-strain elasticity benchmarks. The numerical results show that the method improves substantially over linear spectral approximations and achieves accuracy comparable to snapshot-driven nonlinear manifold constructions.

The remainder of the paper is organized as follows.
Section~\ref{sec:problem_definition} introduces the abstract variational problem, its algebraic discretization, the residual-based error estimate, and the nonlinear manifold reduced-order approximation. Section~\ref{sec:numerical_optimization} presents the residual-driven training strategy, including the alternating optimization algorithm, the affine offline--online decomposition, and the online nonlinear solve. Section~\ref{sec:numerical_experiments} reports numerical results for the advection--diffusion problem and the compressed-plate elasticity benchmark. Finally, Section~\ref{sec:conclusions} summarizes the main findings and discusses open directions.

\section{Problem Definition}
\label{sec:problem_definition}

For a given parameter value $\boldsymbol{\mu}\in\mathcal{D}\subset\mathbb{R}^{S}$, $S\in\mathbb{N}$, consider the following variational problem: find $u(\boldsymbol{\mu})\in\mathcal{V}$ such that
\begin{equation}
a\bigl(u(\boldsymbol{\mu}), v; \boldsymbol{\mu}\bigr) = f(v; \boldsymbol{\mu}), \qquad \forall\, v \in \mathcal{V},
\label{eq:variational_form}
\end{equation}
where $\mathcal{V}\subset H^{1}(\Omega)$ is a Hilbert space defined on the bounded Lipschitz domain $\Omega \subset \mathbb{R}^d$, with $d\in\{1,2,3\}$, endowed with inner product $(\cdot,\cdot)_{\mathcal{V}}$ and associated norm $\|\cdot\|_{\mathcal{V}}$. The mapping $a : \mathcal{V} \times \mathcal{V} \times \mathcal{D} \to \mathbb{R}$ denotes a parametric bilinear form and $f : \mathcal{V} \times \mathcal{D} \to \mathbb{R}$ a parametric linear form. 

For every $\boldsymbol{\mu}\in\mathcal{D}$, we assume that $a(\cdot,\cdot;\boldsymbol{\mu})$ is continuous and coercive. This means that there exist constants $\gamma_a(\boldsymbol{\mu})\le \overline{\gamma}_a<\infty$ and $\alpha_a(\boldsymbol{\mu})\ge \underline{\alpha}_a>0$ such that
\begin{align}
|a(u,v;\boldsymbol{\mu})| &\le \gamma_a(\boldsymbol{\mu})\,\|u\|_{\mathcal{V}}\,\|v\|_{\mathcal{V}} \le \overline{\gamma}_a\,\|u\|_{\mathcal{V}}\,\|v\|_{\mathcal{V}}, && \forall\,u,v\in\mathcal{V}, \label{eq:continuity}\\[4pt]
a(v,v;\boldsymbol{\mu}) &\ge \alpha_a(\boldsymbol{\mu})\,\|v\|_{\mathcal{V}}^{2} \ge \underline{\alpha}_a\,\|v\|_{\mathcal{V}}^{2}, && \forall\,v\in\mathcal{V}. \label{eq:coercivity}
\end{align}
Moreover, for every $\boldsymbol{\mu}\in\mathcal{D}$, the linear form $f(\cdot;\boldsymbol{\mu})$ is bounded: there exists $\gamma_f(\boldsymbol{\mu})\le \overline{\gamma}_f<\infty$ such that
\begin{equation}
|f(v;\boldsymbol{\mu})| \le \gamma_f(\boldsymbol{\mu})\,\|v\|_{\mathcal{V}} \le \overline{\gamma}_f\,\|v\|_{\mathcal{V}}, \qquad \forall\,v\in\mathcal{V}.
\end{equation}
Under these assumptions, the Lax--Milgram theorem~\cite{Brezis2010} guarantees that, for each $\boldsymbol{\mu}\in\mathcal{D}$, variational problem~\eqref{eq:variational_form} admits a unique solution $u(\boldsymbol{\mu})\in\mathcal{V}$.

Given an approximate solution $u_\varepsilon(\boldsymbol{\mu})\in\mathcal{V}$ to~\eqref{eq:variational_form}, we define the associated residual functional $r(\,\cdot\,;u_\varepsilon(\boldsymbol{\mu}),\boldsymbol{\mu}):\mathcal{V}\to\mathbb{R}$ by
\begin{equation}
r(v;u_\varepsilon(\boldsymbol{\mu}),\boldsymbol{\mu}) \coloneq f(v; \boldsymbol{\mu}) - a\bigl(u_\varepsilon(\boldsymbol{\mu}), v; \boldsymbol{\mu}\bigr), \qquad \forall\, v \in \mathcal{V}.
\label{eq:residual}
\end{equation}
Let $e(\boldsymbol{\mu}) \coloneq u(\boldsymbol{\mu}) - u_\varepsilon(\boldsymbol{\mu})$ denote the corresponding error. Then $e(\boldsymbol{\mu})$ satisfies $a(e(\boldsymbol{\mu}),  \allowbreak v; \boldsymbol{\mu}) = r(v;u_\varepsilon(\boldsymbol{\mu}),\boldsymbol{\mu})$ for all $v \in \mathcal{V}$, and the standard residual-based \emph{a posteriori} error estimate reads \cite{rozza2008, veroy2002} 
\begin{equation}
\| e(\boldsymbol{\mu}) \|_{\mathcal{V}} \le \frac{1}{\alpha_a(\boldsymbol{\mu})}\, \| r(\cdot;u_\varepsilon(\boldsymbol{\mu}),\boldsymbol{\mu}) \|_{\mathcal{V}'}.
\label{eq:residual_bound}
\end{equation}

The error is governed by the residual in $\mathcal{V}'$ and is inversely proportional to the coercivity constant, which provides the theoretical framework for the residual-driven estimators and adaptive, snapshot-free strategies described in the subsequent sections.

\subsection{Algebraic Discretization}
\label{sec:matrix_vector}

We now introduce the standard algebraic formulation associated with variational problem~\eqref{eq:variational_form}~\cite{brenner2008, ern2004}. Let $\mathcal{V}_h\subset \mathcal{V}$ be a finite-dimensional subspace with $\dim(\mathcal{V}_h)=N_h\in\mathbb{N}$, and let $\{\chi_i\}_{i=1}^{N_h}\subset \mathcal{V}_h$ be a (hierarchical) basis that is orthonormal with respect to $(\cdot,\cdot)_{\mathcal{V}}$, i.e., $(\chi_i,\chi_j)_{\mathcal{V}}=\delta_{ij}$, where $\delta_{ij}$ denotes the Kronecker delta. We seek the Galerkin approximation $u_h(\boldsymbol{\mu})\in\mathcal{V}_h$ satisfying
\begin{equation}
a\bigl(u_h(\boldsymbol{\mu}), v_h; \boldsymbol{\mu}\bigr) = f(v_h; \boldsymbol{\mu}), \qquad \forall\, v_h \in \mathcal{V}_h.
\label{eq:variational_form_discrete}
\end{equation}
Expanding $u_h(\boldsymbol{\mu})$ in the chosen basis,
\begin{equation}
u_h(\boldsymbol{\mu})=\sum_{j=1}^{N_h} c_j(\boldsymbol{\mu})\,\chi_j, \qquad \mathbf{c}(\boldsymbol{\mu}) \coloneq \bigl[c_1(\boldsymbol{\mu}),\ldots,c_{N_h}(\boldsymbol{\mu})\bigr]^{\mathrm{T}}\in\mathbb{R}^{N_h},
\label{eq:uh_expansion}
\end{equation}
and testing~\eqref{eq:variational_form_discrete} with $\chi_i$, $i=1,\ldots,N_h$, yields the linear system
\begin{equation}
\mathbb{A}_{\boldsymbol{\mu}}\,\mathbf{c}(\boldsymbol{\mu})=\mathbf{f}_{\boldsymbol{\mu}},
\label{eq:matrix_system}
\end{equation}
where the parameter-dependent stiffness matrix $\mathbb{A}_{\boldsymbol{\mu}}\in\mathbb{R}^{N_h\times N_h}$ and load vector $\mathbf{f}_{\boldsymbol{\mu}}\in\mathbb{R}^{N_h}$ are defined entry-wise, for all $i,j\in\{1,\ldots,N_h\}$, by
\begin{equation}
(\mathbb{A}_{\boldsymbol{\mu}})_{ij} \coloneq a(\chi_j,\chi_i;\boldsymbol{\mu}), \qquad (\mathbf{f}_{\boldsymbol{\mu}})_i \coloneq f(\chi_i;\boldsymbol{\mu}).
\label{eq:discrete_A_f}
\end{equation}
The well-posedness of the discrete system in~\eqref{eq:matrix_system} is anchored in the continuity and coercivity of the bilinear form $a(\cdot,\cdot;\boldsymbol{\mu})$, which together ensure the invertibility of $\mathbb{A}_{\boldsymbol{\mu}}$ across all parameters $\boldsymbol{\mu}$ in the parameter domain $\mathcal{D}$.

Given an approximate coefficient vector $\mathbf{c}_\varepsilon(\boldsymbol{\mu})\in\mathbb{R}^{N_h}$ (corresponding to $u_{\varepsilon}(\boldsymbol{\mu})=\sum_{j=1}^{N_h} c_{\varepsilon,j}(\boldsymbol{\mu})\,\chi_j$), we define the discrete residual vector as
\begin{equation}
\mathbf{r}(\boldsymbol{\mu}) \coloneq \mathbf{f}_{\boldsymbol{\mu}} - \mathbb{A}_{\boldsymbol{\mu}}\,\mathbf{c}_\varepsilon(\boldsymbol{\mu}).
\label{eq:residual_vector}
\end{equation}
Its components satisfy $\mathbf{r}_i(\boldsymbol{\mu}) = r(\chi_i;u_\varepsilon(\boldsymbol{\mu}),\boldsymbol{\mu})$, for all $i \in \{1,\ldots,N_h\}$, where $r(\,\cdot\,; \allowbreak u_\varepsilon(\boldsymbol{\mu}),\boldsymbol{\mu})$ denotes the residual functional in~\eqref{eq:residual}.

To relate functional and algebraic norms, we introduce the Gram (mass) matrix $\mathbb{M}\in\mathbb{R}^{N_h\times N_h}$ associated with $(\,\cdot\,,\,\cdot\,)_{\mathcal{V}}$,
\begin{equation}
(\mathbb{M})_{ij} \coloneq (\chi_j,\chi_i)_{\mathcal{V}}.
\label{eq:gram_matrix}
\end{equation}
Under the orthonormality assumption, $\mathbb{M}=\mathbb{I}_{N_h}$, so that $\|\mathbf{x}\|_{\mathbb M}=\|\mathbf{x}\|_2$ and $\|\mathbf{y}\|_{\mathbb M^{-1}}=\|\mathbf{y}\|_2$.

With these definitions, the residual-based \emph{a posteriori} estimate~\eqref{eq:residual_bound} admits the purely algebraic counterpart
\begin{equation}
\|\mathbf{c}(\boldsymbol{\mu})-\mathbf{c}_\varepsilon(\boldsymbol{\mu})\|_{\mathbb{M}} \le \frac{1}{\alpha_a(\boldsymbol{\mu})}\, \|\mathbf{f}_{\boldsymbol{\mu}}-\mathbb{A}_{\boldsymbol{\mu}}\mathbf{c}_\varepsilon(\boldsymbol{\mu})\|_{\mathbb{M}^{-1}}.
\label{eq:algebraic_error_bound}
\end{equation}
This matrix-vector formulation provides the computational realization of the variational problem and serves as the starting point for the model reduction techniques developed in the following sections.

\subsection{Nonlinear Dimensionality Reduction}
\label{sec:nonlinear_dimred}

We next consider a nonlinear representation of the approximate solution $u_\varepsilon(\boldsymbol{\mu})$ at the algebraic level. In particular, we assume that its coefficient vector $\mathbf{c}_\varepsilon(\boldsymbol{\mu})\in\mathbb{R}^{N_h}$ admits a low-dimensional parametrization in terms of latent coordinates $\mathbf{c}_\nu (\boldsymbol{\mu}) \in \mathbb{R}^{N_\nu}$, with $N_\nu \ll N_h$, via the nonlinear manifold ansatz
\begin{equation}
\label{eq:nonlinear_manifold}
\mathbf{c}_\varepsilon(\boldsymbol{\mu}) = \mathbb{V}\,\mathbf{c}_\nu(\boldsymbol{\mu}) + \mathbb{U}\,\mathbb{X}\,\mathbf{h}\!\left(\mathbf{c}_\nu(\boldsymbol{\mu})\right).
\end{equation}
Here, $\mathbb{X}\in\mathbb{R}^{(N_\varepsilon-N_\nu)\times N_\ell}$ denotes the matrix of lifting coefficients to be inferred numerically, where $N_\varepsilon\in\mathbb{N}$ is the dimension of the lifting space. The matrices $\mathbb{V}\in\mathbb{R}^{N_h\times N_\nu}$ and $\mathbb{U}\in\mathbb{R}^{N_h\times (N_\varepsilon-N_\nu)}$ have $\mathbb{M}$-orthonormal columns and span complementary subspaces of $\mathbb{R}^{N_h}$. Finally, $\mathbf{h}:\mathbb{R}^{N_\nu}\rightarrow\mathbb{R}^{N_\ell}$ is a prescribed nonlinear feature map of the latent coordinates, typically encoding higher-order, e.g., polynomial, interactions, where $N_\ell$ denotes the number of nonlinear features. A specific quadratic choice of $\mathbf{h}$, with the corresponding value of $N_\ell$, is given in~\eqref{eq:quadratic_feature_map}.

Substituting~\eqref{eq:nonlinear_manifold} into the full-order system~\eqref{eq:matrix_system} and projecting onto the latent subspace spanned by the columns of $\mathbb{V}$ yields the nonlinear reduced-order model
\begin{equation}
\label{eq:nonlinear_rom}
\mathbb{A}_{\boldsymbol{\mu}}^\dagger\,\mathbf{c}_\nu(\boldsymbol{\mu}) + \mathbb{A}_{\boldsymbol{\mu}}^\ddagger \, \mathbb{X} \, \mathbf{h}\!\left(\mathbf{c}_\nu(\boldsymbol{\mu})\right) = \mathbf{f}_{\boldsymbol{\mu}}^\dagger \qquad \forall \boldsymbol{\mu}\in\mathcal{D},
\end{equation}
where the projected operators are defined by
\begin{equation}
\label{eq:projected_operators}
\mathbb{A}_{\boldsymbol{\mu}}^\dagger \coloneq \mathbb{V}^\top \mathbb{A}_{\boldsymbol{\mu}} \mathbb{V}, \qquad \mathbb{A}_{\boldsymbol{\mu}}^\ddagger \coloneq \mathbb{V}^\top \mathbb{A}_{\boldsymbol{\mu}} \mathbb{U}, \qquad \mathbf{f}_{\boldsymbol{\mu}}^\dagger \coloneq \mathbb{V}^\top \mathbf{f}_{\boldsymbol{\mu}}.
\end{equation}

Recall that $\mathbb{V}$ and $\mathbb{U}$ are $\mathbb{M}$-orthonormal, i.e., 
\begin{equation*}
\mathbb{V}^\top \mathbb{M} \mathbb{V} = \mathbb{I}_{N_\nu}, \qquad \mathbb{U}^\top \mathbb{M} \mathbb{U} = \mathbb{I}_{N_\varepsilon-N_\nu}, \qquad \mathbb{U}^\top \mathbb{M} \mathbb{V} = \mathbf{0}.
\end{equation*}
In principle, these two matrices may be jointly optimized to obtain an accurate manifold parametrization. In the present study, however, a simplified construction is adopted. After expressing the solution in the chosen $\mathbb{M}$-orthonormal hierarchical basis, the columns of $\mathbb{V}$ and $\mathbb{U}$ are selected as coordinate vectors in that modal basis, that is, $\mathbb{I}_{N_h} = [\, \mathbb{V}, \, \mathbb{U}, \, \ldots \,]$. Thus, $\mathbb{V}$ corresponds to the first $N_\nu$ modes, while $\mathbb{U}$ contains the subsequent $N_\varepsilon-N_\nu$ modes. This construction provides a baseline for more advanced enrichment strategies, such as the greedy approach for quadratic manifolds proposed in~\cite{schwerdtner2024greedy}. While that method adaptively builds the basis through greedy updates, the present study focuses on evaluating the residual-driven training framework within a prescribed hierarchical subspace.

The matrix $\mathbb{X}$ specifies the lifting coefficients and governs the manifold's curvature through~\eqref{eq:nonlinear_manifold}. Existing data-driven approaches \cite{Barnett2022, Geelen2024} estimate $\mathbb{X}$ from a training set of high-fidelity solutions of~\eqref{eq:matrix_system}. In contrast, in this work we identify $\mathbb{X}$ using only reduced-dimension operations on precomputed affine blocks, without any full-order solves.

\subsection{Optimal Lifting Matrix}
\label{subsec:optimal_lifting_matrix}

We define the optimal lifting matrix $\mathbb{X}$ as the one that yields, on average over the parameter domain $\mathcal{D}$, the most accurate approximation of the full-order coefficient vector $\mathbf{c}(\boldsymbol{\mu})$ by the manifold representation~\eqref{eq:nonlinear_manifold}. Assuming a uniform weighting on $\mathcal{D}$, this leads to a constrained optimization problem of the form:
\begin{equation}
\label{eq:exact_optimization}
\begin{aligned}
\mathbb{X} &= \operatorname*{arg\,min}_{\tilde{\mathbb{X}} \in \mathbb{R}^{(N_\varepsilon - N_\nu) \times N_\ell}} \int_{\mathcal{D}} \big\| \mathbf{c}(\boldsymbol{\mu}) - \mathbb{V} \, \mathbf{c}_\nu (\boldsymbol{\mu}) - \mathbb{U} \, \tilde{\mathbb{X}} \, \mathbf{h} \! \left(\mathbf{c}_\nu (\boldsymbol{\mu}) \right) \big\|_{\mathbb{M}}^{2}\, d\boldsymbol{\mu} \\
&\text{s.t.}\quad \mathbb{A}^\dagger_{\boldsymbol{\mu}}\,\mathbf{c}_\nu(\boldsymbol{\mu}) + \mathbb{A}^\ddagger_{\boldsymbol{\mu}} \, \tilde{\mathbb{X}} \, \mathbf{h} \! \left( \mathbf{c}_\nu(\boldsymbol{\mu}) \right) = \mathbf{f}^\dagger_{\boldsymbol{\mu}}, \qquad \forall \, \boldsymbol{\mu} \in \mathcal{D}.
\end{aligned}
\end{equation}
The constraint ensures that the latent coordinates used to fit the manifold are not arbitrary encoding variables, but are the coordinates produced by the nonlinear reduced Galerkin system that will be solved online.

If full-order solution snapshots $\mathbf{c}(\boldsymbol{\mu})$ are available, Problem~\eqref{eq:exact_optimization} can be approximated in the offline stage by replacing the parameter-domain integral with a finite training average and minimizing the resulting empirical objective. In this work, however, the goal is to construct a snapshot-free alternative that avoids computing $\mathbf{c}(\boldsymbol{\mu})$ over the parameter domain $\mathcal{D}$.

Using $\mathbf{r}(\boldsymbol{\mu}) =
\mathbf{f}_{\boldsymbol{\mu}} -\mathbb{A}_{\boldsymbol{\mu}}\mathbf{c}_\varepsilon(\boldsymbol{\mu})$, the residual-based estimate~\eqref{eq:algebraic_error_bound} can be integrated over $\mathcal{D}$ to obtain
\begin{equation}
\label{eq:bounded_cost_function}
\int_{\mathcal{D}} \| \mathbf{c}(\boldsymbol{\mu}) - \mathbf{c}_\varepsilon(\boldsymbol{\mu}) \|_{\mathbb{M}}^2 \, d\boldsymbol{\mu} \le \int_{\mathcal{D}} \frac{1}{\alpha_a(\boldsymbol{\mu})^2} \, \| \mathbf{f}_{\boldsymbol{\mu}} - \mathbb{A}_{\boldsymbol{\mu}} \mathbf{c}_\varepsilon(\boldsymbol{\mu}) \|_{\mathbb{M}^{-1}}^2 \, d\boldsymbol{\mu}.
\end{equation}
The estimate~\eqref{eq:bounded_cost_function} provides a computable surrogate for the inaccessible average state error. Therefore, instead of minimizing the true error, which requires the full-order solution $\mathbf{c}(\boldsymbol{\mu})$, the lifting matrix is identified by minimizing the residual-based upper bound on the right-hand side of~\eqref{eq:bounded_cost_function}. Substituting the nonlinear manifold ansatz~\eqref{eq:nonlinear_manifold} yields the following snapshot-free optimization problem:
\begin{equation}
\label{eq:bound_optimization}
\begin{aligned}
\mathbb{X} &= \operatorname*{arg\,min}_{\tilde{\mathbb{X}} \in \mathbb{R}^{(N_\varepsilon - N_\nu) \times N_\ell}} \int_{\mathcal{D}} \frac{1}{\alpha_a(\boldsymbol{\mu})^2} \, \big\| \mathbf{f}_{\boldsymbol{\mu}} - \mathbb{A}_{\boldsymbol{\mu}} \big( \mathbb{V}\,\mathbf{c}_\nu(\boldsymbol{\mu}) + \mathbb{U}\,\tilde{\mathbb{X}} \, \mathbf{h} \! \left(\mathbf{c}_\nu (\boldsymbol{\mu}) \right) \big) \big\|_{\mathbb{M}^{-1}}^{2} \, d\boldsymbol{\mu} \\
&\text{s.t.}\quad \mathbb{A}^\dagger_{\boldsymbol{\mu}}\,\mathbf{c}_\nu(\boldsymbol{\mu}) + \mathbb{A}^\ddagger_{\boldsymbol{\mu}} \, \tilde{\mathbb{X}} \, \mathbf{h} \! \left( \mathbf{c}_\nu (\boldsymbol{\mu}) \right) = \mathbf{f}^\dagger_{\boldsymbol{\mu}}, \qquad \forall\,\boldsymbol{\mu}\in\mathcal{D}.
\end{aligned}
\end{equation}
In this formulation, the lifting matrix is determined by minimizing a computable residual-based surrogate for the average state error, without requiring full-order solution snapshots.

\section{Numerical Optimization}
\label{sec:numerical_optimization}

Problem~\eqref{eq:bound_optimization} is a constrained minimization: the lifting matrix $\mathbb{X}$ enters both the residual objective and, through the reduced constraint, the latent coordinates $\mathbf{c}_\nu(\boldsymbol{\mu})$. We address it with a fixed-point procedure for the associated first-order residual-training equations. In the spirit of block-coordinate methods~\cite{Tseng2001}, the procedure alternates between (i) updating the latent coordinates and (ii) updating the lifting matrix; the complete offline training procedure is summarized in Algorithm~\ref{alg:nlmanifold_training_online}.

We begin by introducing a discrete counterpart of~\eqref{eq:bound_optimization}, obtained by quadrature over the parameter domain. Let $\{(\boldsymbol{\mu}_k,\omega_k)\}_{k=1}^{N_\omega}$ be a set of parameter samples and associated positive quadrature weights, with $N_\omega\in\mathbb{N}$. For instance, a Monte Carlo rule draws $\boldsymbol{\mu}_k$ uniformly from $\mathcal{D}$ and sets $\omega_k = 1/N_\omega$. Approximating the integral in~\eqref{eq:bound_optimization} by this quadrature rule yields the problem
\begin{equation}
\label{eq:bound_optimization_discrete}
\begin{aligned}
\vspace{-0.1cm}
\mathbb{X} &= \operatorname*{arg\,min}_{\tilde{\mathbb{X}} \in \mathbb{R}^{(N_\varepsilon - N_\nu)\times N_\ell}} \sum_{k=1}^{N_\omega} \frac{\omega_k}{\alpha_a(\boldsymbol{\mu}_k)^2} \big\| \mathbf{f}_{\boldsymbol{\mu}_k} - \mathbb{A}_{\boldsymbol{\mu}_k} \big( \mathbb{V}\,\mathbf{c}_\nu(\boldsymbol{\mu}_k) + \mathbb{U} \, \tilde{\mathbb{X}} \, \mathbf{h} \! \left(\mathbf{c}_\nu (\boldsymbol{\mu}_k) \right) \big) \big\|_{\mathbb{M}^{-1}}^2 \\
&\text{s.t.}\quad \mathbb{A}^\dagger_{\boldsymbol{\mu}_k} \, \mathbf{c}_\nu (\boldsymbol{\mu}_k) + \mathbb{A}^\ddagger_{\boldsymbol{\mu}_k} \, \tilde{\mathbb{X}} \, \mathbf{h} \! \left( \mathbf{c}_\nu(\boldsymbol{\mu}_k) \right) = \mathbf{f}^\dagger_{\boldsymbol{\mu}_k}, \qquad \forall k\in \{1, \ldots, N_\omega \}.
\end{aligned}
\end{equation}

Assuming the iteration below admits a fixed point (analysis is deferred to future work), we initialize $\mathbb{X}_0=\mathbf{0}$ and alternate between the following two steps until convergence is reached:
\begin{enumerate}[label=(\roman*), leftmargin=*, itemsep=0.3em, topsep=0.3em]
\item \textbf{Latent-space update.}
For each sample $\boldsymbol{\mu}_k$, compute the latent coordinates $\mathbf{c}_{\nu,n}(\boldsymbol{\mu}_k)\in\mathbb{R}^{N_\nu}$ as the solution of the nonlinear reduced system
\begin{equation}
\label{eq:latent_update}
\mathbb{A}^\dagger_{\boldsymbol{\mu}_k}\,\mathbf{c}_{\nu,n}(\boldsymbol{\mu}_k) + \mathbb{A}^\ddagger_{\boldsymbol{\mu}_k} \, \mathbb{X}_{n} \, \mathbf{h} \! \left(\mathbf{c}_{\nu,n} (\boldsymbol{\mu}_k) \right) = \mathbf{f}^\dagger_{\boldsymbol{\mu}_k},
\end{equation}
\item \textbf{Lifting-matrix update.}
With $\{\mathbf{c}_{\nu,n}(\boldsymbol{\mu}_k)\}_{k=1}^{N_\omega}$ fixed, update $\mathbb{X}_{n+1}$ by solving the weighted least-squares problem
\begin{equation}
\label{eq:lifting_update}
\hspace{-0.6cm}
\vspace{-0.1cm}
\begin{aligned}
\mathbb{X}_{n+1} = \operatorname*{arg\,min}_{\tilde{\mathbb{X}} \in \mathbb{R}^{(N_\varepsilon - N_\nu)\times N_\ell}} \sum_{k=1}^{N_\omega} \frac{\omega_k}{\alpha_a(\boldsymbol{\mu}_k)^2} \Big\| \mathbf{f}_{\boldsymbol{\mu}_k} - \mathbb{A}_{\boldsymbol{\mu}_k} \! \big( \mathbb{V}\,\mathbf{c}_{\nu,n}(\boldsymbol{\mu}_k) + \! \mathbb{U}\,\tilde{\mathbb{X}}\,\mathbf{h}\!\left(\mathbf{c}_{\nu,n}(\boldsymbol{\mu}_k)\right)  \big)\Big\|_{\mathbb{M}^{-1}}^2.
\end{aligned}
\end{equation}
\end{enumerate}
where $n \in \mathbb{N}$ is the iteration counter.

Convergence is monitored through the change in the reconstructed reduced coefficients induced by the lifting update, with the current latent coordinates held fixed. By $\mathbb{M}$-orthogonality of $\mathbb{V}$ and $\mathbb{U}$, this criterion reduces to the low-dimensional condition
\begin{equation}
\label{eq:conv_criterion}
\max_{1\leq k \leq N_\omega} \frac{\|\mathbb{U}(\mathbb{X}_{n+1}-\mathbb{X}_{n})\mathbf{h}\left(\mathbf{c}_{\nu,n}(\boldsymbol{\mu}_k)\right)\|^{2}_{\mathbb{M}}}{\|\mathbb{V}\mathbf{c}_{\nu,n}(\boldsymbol{\mu}_k)\|^{2}_{\mathbb{M}} + \|\mathbb{U}\mathbb{X}_{n+1}\mathbf{h}\left(\mathbf{c}_{\nu,n}(\boldsymbol{\mu}_k)\right)\|^{2}_{\mathbb{M}}} \le \mathrm{tol}^{2},
\end{equation}
which is the form used in Algorithm~\ref{alg:nlmanifold_training_online}. The mean weighted residual norm associated with~\eqref{eq:bound_optimization_discrete} is also tracked during training to verify the decrease of the objective.

\subsection{Affine Parametric Dependence}
\label{sec:affine_dependence}

We now assume that the operators in~\eqref{eq:variational_form}, and hence the stiffness matrix and load vector in~\eqref{eq:matrix_system}, admit an \emph{affine} parametric decomposition, as commonly required to enable efficient offline--online evaluations in projection-based model reduction~\cite{barrault2004, grepl2007, quarteroni2015}. In particular, for every $\boldsymbol{\mu}\in\mathcal{D}$, we write
\begin{equation}
\label{eq:affine_decomposition}
\mathbb{A}_{\boldsymbol{\mu}} = \sum_{p=1}^{N_A} \Theta_p^a(\boldsymbol{\mu}) \,\mathbb{A}_{p}, \qquad \mathbf{f}_{\boldsymbol{\mu}} = \sum_{q=1}^{N_F} \Theta_q^f(\boldsymbol{\mu}) \,\mathbf{f}_q,
\end{equation}
where the coefficient functions $\Theta_p^a,\Theta_q^f:\mathcal{D}\to\mathbb{R}$ are known and the components $\mathbb{A}_p\in\mathbb{R}^{N_h\times N_h}$ and $\mathbf{f}_q\in\mathbb{R}^{N_h}$ are parameter-independent for all $p\in\{1,\ldots,N_A\}$ and $q\in\{1,\ldots,N_F\}$.

This structure allows offline precomputation of the reduced operators entering \eqref{eq:nonlinear_rom} and in the objective in~\eqref{eq:bound_optimization_discrete}. In particular, using~\eqref{eq:affine_decomposition} and~\eqref{eq:projected_operators}, we obtain the affine reduced decompositions
\begin{equation}
\label{eq:reduced_affine}
\mathbb{A}^\dagger_{\boldsymbol{\mu}} = \sum_{p=1}^{N_A} \Theta^a_p(\boldsymbol{\mu}) \, \mathbb{A}_p^{\dagger}, \qquad \mathbb{A}_{\boldsymbol{\mu}}^\ddagger = \sum_{p=1}^{N_A} \Theta^a_p(\boldsymbol{\mu}) \, \mathbb{A}_p^{\ddagger}, \qquad \mathbf{f}^\dagger_{\boldsymbol{\mu}} = \sum_{q=1}^{N_F} \Theta^f_q(\boldsymbol{\mu}) \, \mathbf{f}_q^{\dagger},
\end{equation}
with parameter-independent reduced blocks
\begin{equation}
\label{eq:reduced_blocks}
\mathbb{A}_p^{\dagger} \coloneq \mathbb{V}^\top \mathbb{A}_p \mathbb{V}, \qquad \mathbb{A}_p^{\ddagger} \coloneq \mathbb{V}^\top \mathbb{A}_p \mathbb{U}, \qquad \mathbf{f}_q^{\dagger} \coloneq \mathbb{V}^\top \mathbf{f}_q.
\end{equation}
To enable an efficient offline--online treatment of the objective in~\eqref{eq:bound_optimization_discrete}, consider a single quadrature sample $\boldsymbol{\mu}_k$ and introduce the abbreviations $\mathbf{c}_{\nu,\boldsymbol{\mu}_k} \coloneq \mathbf{c}_\nu(\boldsymbol{\mu}_k)$ and $\mathbf{h}_{\boldsymbol{\mu}_k} \coloneq \mathbf{h}(\mathbf{c}_{\nu,\boldsymbol{\mu}_k})$. Recalling~\eqref{eq:nonlinear_manifold}, $\mathbf{c}_\varepsilon(\boldsymbol{\mu}_k) = \mathbb{V}\mathbf{c}_{\nu,\boldsymbol{\mu}_k} + \mathbb{U}\mathbb{X} \mathbf{h}_{\boldsymbol{\mu}_k}$, we introduce the latent-space residual
\begin{equation*}
\mathbf{b}_{\boldsymbol{\mu}_k} \coloneq \mathbf{f}_{\boldsymbol{\mu}_k} - \mathbb{A}_{\boldsymbol{\mu}_k}\,\mathbb{V}\,\mathbf{c}_{\nu,\boldsymbol{\mu}_k},   
\end{equation*}
which corresponds to the residual obtained from the latent contribution alone. The full residual associated with $\mathbb{X}$ can then be written as $\mathbf{r}(\boldsymbol{\mu}_k;\mathbb{X}) = \mathbf{b}_{\boldsymbol{\mu}_k} - \mathbb{A}_{\boldsymbol{\mu}_k}\, \mathbb{U} \,\mathbb{X} \,\mathbf{h}_{\boldsymbol{\mu}_k}$. Expanding its weighted norm yields the quadratic form
\begin{equation}
\label{eq:quad_expand_master}
\begin{aligned}
\big\|\mathbf{r}(\boldsymbol{\mu}_k;\mathbb{X})\big\|_{\mathbb{M}^{-1}}^{2} \coloneq \gamma_{\boldsymbol{\mu}_k} - 2\,\mathbf{h}_{\boldsymbol{\mu}_k}^\top \mathbb{X}^\top \mathbf{d}_{\boldsymbol{\mu}_k} + \mathbf{h}_{\boldsymbol{\mu}_k}^\top \mathbb{X}^\top \mathbb{E}_{\boldsymbol{\mu}_k} \mathbb{X}\,\mathbf{h}_{\boldsymbol{\mu}_k} ,
\end{aligned}
\end{equation}
where 
\begin{equation*}
  \gamma_{\boldsymbol{\mu}_k} \coloneq \mathbf{b}_{\boldsymbol{\mu}_k}^\top \mathbb{M}^{-1}\mathbf{b}_{\boldsymbol{\mu}_k}, \quad\; \mathbf{d}_{\boldsymbol{\mu}_k} \coloneq \mathbb{U}^\top \mathbb{A}_{\boldsymbol{\mu}_k}^\top \mathbb{M}^{-1}\mathbf{b}_{\boldsymbol{\mu}_k}, \quad\; \mathbb{E}_{\boldsymbol{\mu}_k} \coloneq \mathbb{U}^\top \mathbb{A}_{\boldsymbol{\mu}_k}^\top \mathbb{M}^{-1}\mathbb{A}_{\boldsymbol{\mu}_k}\mathbb{U}.
\end{equation*}
Since $\gamma_{\boldsymbol{\mu}_k}$ is independent of $\mathbb{X}$, it can be omitted from the minimization.

Using the affine decomposition~\eqref{eq:affine_decomposition}, the vectors $\mathbf{d}_{\boldsymbol{\mu}_k}$ and matrices $\mathbb{E}_{\boldsymbol{\mu}_k}$ in~\eqref{eq:quad_expand_master} can be assembled from parameter-independent blocks. Define, for $i,j\in\{1,\ldots,N_A\}$ and $\ell\in\{1,\ldots,N_F\}$,
\begin{equation}
\label{eq:offline_blocks}
\begin{aligned}
\mathbb{S}_{ij} &\coloneq \mathbb{U}^\top \mathbb{A}_i^\top \mathbb{M}^{-1} \mathbb{A}_j \mathbb{U} \in \mathbb{R}^{(N_\varepsilon - N_\nu)\times(N_\varepsilon - N_\nu)},\\
\mathbb{P}_{ij} &\coloneq \mathbb{U}^\top \mathbb{A}_i^\top \mathbb{M}^{-1} \mathbb{A}_j \mathbb{V} \in \mathbb{R}^{(N_\varepsilon - N_\nu)\times N_\nu},\\
\mathbf{q}_{i\ell} &\coloneq \mathbb{U}^\top \mathbb{A}_i^\top \mathbb{M}^{-1} \mathbf{f}_\ell \in \mathbb{R}^{N_\varepsilon - N_\nu}.
\end{aligned}
\end{equation}
Then, for each quadrature sample $\boldsymbol{\mu}_k$, the coefficients in~\eqref{eq:quad_expand_master} admit the affine expansions
\begin{align}
\label{eq:E_d_affine}
\mathbb{E}_{\boldsymbol{\mu}_k} &= \sum_{i=1}^{N_A}\sum_{j=1}^{N_A} \Theta_i^a(\boldsymbol{\mu}_k)\,\Theta_j^a(\boldsymbol{\mu}_k)\,\mathbb{S}_{ij},\\
\mathbf{d}_{\boldsymbol{\mu}_k} &= \sum_{i=1}^{N_A}\sum_{\ell=1}^{N_F} \Theta_i^a(\boldsymbol{\mu}_k)\,\Theta_\ell^f(\boldsymbol{\mu}_k)\,\mathbf{q}_{i\ell} \;-\; \sum_{i=1}^{N_A}\sum_{j=1}^{N_A} \Theta_i^a(\boldsymbol{\mu}_k)\,\Theta_j^a(\boldsymbol{\mu}_k)\,\mathbb{P}_{ij}\,\mathbf{c}_{\nu,{\boldsymbol{\mu}_k}}.
\end{align}
All blocks in~\eqref{eq:offline_blocks} are parameter-independent and can be precomputed before the training phase. In the training phase, only the scalar coefficients $\{\Theta_i^a(\boldsymbol{\mu}_k)\}_{i=1}^{N_A}$ and $\{\Theta_\ell^f(\boldsymbol{\mu}_k)\}_{\ell=1}^{N_F}$ are needed to assemble $\mathbb{E}_{\boldsymbol{\mu}_k}$ and $\mathbf{d}_{\boldsymbol{\mu}_k}$, together with the current latent coordinates $\mathbf{c}_{\nu,\boldsymbol{\mu}_k}$. The assembly cost is independent of $N_h$; in particular, forming $\mathbb{E}_{\boldsymbol{\mu}_k}$ scales as $\mathcal{O} \left(N_A^2 (N_\varepsilon-N_\nu)^2\right)$, while forming $\mathbf{d}_{\boldsymbol{\mu}_k}$ scales as $\mathcal{O} ( N_A( N_F + N_A N_\nu )(N_\varepsilon - N_\nu) )$. These operations are repeated for all $k \in \{ 1,\ldots,N_\omega \}$ at each iteration of the alternating scheme.

Substituting~\eqref{eq:quad_expand_master} into~\eqref{eq:lifting_update} (and dropping the constant term $\gamma_{\boldsymbol{\mu}_k}$) shows that the lifting update amounts to minimizing a weighted sum of convex quadratic forms in $\tilde{\mathbb{X}}$:
\begin{equation}
\label{eq:lifting_update_simplified}
\mathbb{X}_{n+1} = \operatorname*{arg\,min}_{\tilde{\mathbb{X}} \in \mathbb{R}^{(N_\varepsilon - N_\nu)\times N_\ell}} \sum_{k=1}^{N_\omega} \frac{\omega_k}{\alpha_a(\boldsymbol{\mu}_k)^2} \left( \mathbf{h}_{\boldsymbol{\mu}_k}^\top \tilde{\mathbb{X}}^\top \mathbb{E}_{\boldsymbol{\mu}_k} \tilde{\mathbb{X}}\,\mathbf{h}_{\boldsymbol{\mu}_k} -2\,\mathbf{h}_{\boldsymbol{\mu}_k}^\top \tilde{\mathbb{X}}^\top \mathbf{d}_{\boldsymbol{\mu}_k}
\right) .
\end{equation}

The coercivity factors in~\eqref{eq:lifting_update_simplified} are weights inherited from the residual-based error estimate, but they are not essential for the snapshot-free identification of $\mathbb{X}$. If the problem is known to be coercive on $\mathcal{D}$ but no sharp or easily computable expression for $\alpha_a(\boldsymbol{\mu})$ is available, one may instead minimize the unweighted residual objective, using only the quadrature weights $\omega_k$. This corresponds to replacing the coercivity factor by the same positive constant for all training samples, which leaves the minimizer unchanged. The resulting training problem remains fully computable, although its objective values no longer define certified \emph{a posteriori} error bounds unless a valid coercivity lower bound is available.

Setting the gradient of~\eqref{eq:lifting_update_simplified} with respect to $\tilde{\mathbb{X}}$ to zero yields the normal equations $\mathbb{G}\,\mathrm{vec}(\mathbb{X}_{n+1}) = \sum_{k=1}^{N_\omega}\varpi_k\,\mathrm{vec}(\mathbf{d}_{\boldsymbol{\mu}_k} \mathbf{h}_{\boldsymbol{\mu}_k}^\top)$, with $\varpi_k\coloneq\omega_k/\alpha_a(\boldsymbol{\mu}_k)^2$ and
\begin{equation}
\label{eq:normal_eqs}
\mathbb{G}\coloneq \sum_{k=1}^{N_\omega} \varpi_k\,(\mathbf{h}_{\boldsymbol{\mu}_k} \mathbf{h}_{\boldsymbol{\mu}_k}^\top) \otimes \mathbb{E}_{\boldsymbol{\mu}_k}.
\end{equation}
In practice, explicitly assembling the Kronecker matrix in~\eqref{eq:normal_eqs} is unnecessary and may be prohibitively expensive. Instead, it is convenient to work with the equivalent matrix-form normal equation
\begin{equation}
\label{eq:doable_equation}
\sum_{k=1}^{N_\omega} \varpi_k \left( \mathbb{E}_{\boldsymbol{\mu}_k} \mathbb{X}_{n+1} \mathbf{h}_{\boldsymbol{\mu}_k} \right)\mathbf{h}_{\boldsymbol{\mu}_k}^\top = \sum_{k=1}^{N_\omega} \varpi_k  \mathbf{d}_{\boldsymbol{\mu}_k}\mathbf{h}_{\boldsymbol{\mu}_k}^\top.
\end{equation}
When the coercivity-independent formulation is adopted, \eqref{eq:doable_equation} is solved with weights $\omega_k$ in place of $\varpi_k$.

Note that the operator $\mathbb{G}$ in~\eqref{eq:normal_eqs} is symmetric positive semidefinite, since $\mathrm{vec}(\mathbb{Z})^\top \mathbb{G}\,\mathrm{vec}(\mathbb{Z}) = \sum_{k=1}^{N_\omega} \varpi_k \, \| \mathbb{M}^{-1/2} \mathbb{A}_{\boldsymbol{\mu}_k} \mathbb{U} \, \mathbb{Z} \, \mathbf{h}_{\boldsymbol{\mu}_k} \|_2^2 \ge 0$ for all $\mathbb{Z} \in \mathbb{R}^{(N_\varepsilon-N_\nu) \times N_\ell}$. Coercivity of $a(\cdot,\cdot;\boldsymbol{\mu}_k)$ together with $\mathbb{U}$ having full column rank makes $\mathbb{A}_{\boldsymbol{\mu}_k}\mathbb{U}$ injective, so each summand vanishes if and only if $\mathbb{Z}\,\mathbf{h}_{\boldsymbol{\mu}_k}=\mathbf{0}$. Hence $\mathbb{G}$ is positive definite (and the lifting update is uniquely solvable) if and only if the feature vectors $\smash{\{ \mathbf{h}_{\boldsymbol{\mu}_k} \}_{k=1}^{N_\omega}}$ span $\mathbb{R}^{N_\ell}$, which requires at least $N_\omega \ge N_\ell$ samples. When this fails, a minimum-norm solution can still be recovered by adding a small Tikhonov regularization term.

Because $\mathbb{G}$ is symmetric positive semidefinite, the normal equations can be solved by a conjugate-gradient method. The unknown is stored in matrix form, while the iterative solver is applied to the corresponding vectorized normal equations. At each iteration, the action of $\mathbb{G}$ on a trial matrix $\mathbb{Z}\in\mathbb{R}^{(N_\varepsilon-N_\nu)\times N_\ell}$ is computed by evaluating the left-hand side of~\eqref{eq:doable_equation} with $\mathbb{X}_{n+1}$ replaced by $\mathbb{Z}$. The resulting matrix is then vectorized and used as the matrix--vector product required by the iterative solver. Thus, the full Kronecker matrix is never assembled or stored.

\begin{algorithm}
\caption{Snapshot-free identification of the nonlinear-manifold lifting matrix}
\label{alg:nlmanifold_training_online}
\begin{algorithmic}[1]
\Require Bases $\mathbb{V}$, $\mathbb{U}$, feature map $\mathbf{h}(\cdot)$, quadrature samples and weights $\{(\boldsymbol{\mu}_k, \varpi_k \coloneq \omega_k / \alpha_a(\boldsymbol{\mu}_k)^2)\}_{k=1}^{N_\omega}$, affine components $\{\mathbb{A}_p\}_{p=1}^{N_A}$, $\{\mathbf{f}_q\}_{q=1}^{N_F}$ and coefficient functions $\Theta_p^a(\cdot)$, $\Theta_q^f(\cdot)$, mass matrix $\mathbb{M}$, tolerance $\mathrm{tol}$.
\Ensure Trained lifting matrix $\mathbb{X}$ 
\vspace{0.3em}

%\Statex \textbf{Precomputation}
\State Compute $\mathbb{M}^\dagger = \mathbb{V}^\top \mathbb{M} \mathbb{V}$ \,and\, $\mathbb{M}^\ddagger = \mathbb{U}^\top \mathbb{M} \mathbb{U}$ 
\For{$p=1,\dots,N_A$}
    \State Compute $\mathbb{A}_p^\dagger = \mathbb{V}^\top \mathbb{A}_p \mathbb{V}$ and $\mathbb{A}_p^\ddagger = \mathbb{V}^\top \mathbb{A}_p \mathbb{U}$
\EndFor
\For{$q=1,\dots,N_F$}
    \State Compute $\mathbf{f}_q^\dagger = \mathbb{V}^\top \mathbf{f}_q$
\EndFor
\For{$i=1,\dots,N_A$}
    \For{$j=1,\dots,N_A$}
        \State Compute $\mathbb{S}_{ij} = \mathbb{U}^\top \mathbb{A}_i^\top \mathbb{M}^{-1}\mathbb{A}_j\mathbb{U}$ and  $\mathbb{P}_{ij} = \mathbb{U}^\top \mathbb{A}_i^\top \mathbb{M}^{-1}\mathbb{A}_j\mathbb{V}$
    \EndFor
    \For{$\ell=1,\dots,N_F$}
        \State Compute $\mathbf{q}_{i\ell} = \mathbb{U}^\top \mathbb{A}_i^\top \mathbb{M}^{-1}\mathbf{f}_\ell$
    \EndFor
\EndFor

\vspace{0.2cm}
\State Initialize $\mathbb{X}_0 = \mathbf{0}$
\For{$n=0,1,2,\ldots$}
\vspace{0.2cm}

    \Statex \hspace{12pt} \textbf{Step 1: latent-space update}
    \For{$k=1,\dots,N_\omega$}
        \State Assemble $\mathbb{A}_{\boldsymbol{\mu}_k}^\dagger = \sum_{p=1}^{N_A} \Theta_p^a(\boldsymbol{\mu}_k)\mathbb{A}_p^\dagger$ \,and\, $\mathbb{A}_{\boldsymbol{\mu}_k}^\ddagger = \sum_{p=1}^{N_A} \Theta_p^a(\boldsymbol{\mu}_k)\mathbb{A}_p^\ddagger$
        \State Assemble $ \mathbf{f}_{\boldsymbol{\mu}_k}^\dagger = \sum_{q=1}^{N_F}\Theta_q^f(\boldsymbol{\mu}_k)\mathbf{f}_q^\dagger$
        \State Solve for $\mathbf{c}_{\nu,n}(\boldsymbol{\mu}_k)$: $\mathbb{A}_{\boldsymbol{\mu}_k}^\dagger \mathbf{c}_{\nu,n}(\boldsymbol{\mu}_k) + \mathbb{A}_{\boldsymbol{\mu}_k}^\ddagger \mathbb{X}_n \mathbf{h}\!\left(\mathbf{c}_{\nu,n}(\boldsymbol{\mu}_k)\right) = \mathbf{f}_{\boldsymbol{\mu}_k}^\dagger$
    \EndFor

    \vspace{0.1cm}
    \Statex \hspace{12pt} \textbf{Step 2: lifting-matrix update}
    \For{$k=1,\dots,N_\omega$}
        \State Set $\mathbf{h}_k=\mathbf{h}(\mathbf{c}_{\nu,n}(\boldsymbol{\mu}_k))$
        \State Assemble $\mathbb{E}_{\boldsymbol{\mu}_k} \! = \sum_{i,j=1}^{N_A, N_A}\Theta_i^a(\boldsymbol{\mu}_k)\Theta_j^a(\boldsymbol{\mu}_k)\mathbb{S}_{ij}$
        \State Assemble \mbox{$\mathbf{d}_{\boldsymbol{\mu}_k} \!= \sum_{i,\ell=1}^{N_A, N_F} \Theta_i^a(\boldsymbol{\mu}_k)\Theta_\ell^f(\boldsymbol{\mu}_k)\mathbf{q}_{i\ell} - \sum_{i,j=1}^{N_A, N_A} \Theta_i^a(\boldsymbol{\mu}_k)\Theta_j^a(\boldsymbol{\mu}_k) \mathbb{P}_{ij}\mathbf{c}_{\nu,n}(\boldsymbol{\mu}_k)$}
    \EndFor
    \State Solve for $\mathbb{X}_{n+1}$: $\sum_{k=1}^{N_\omega} \varpi_k \left( \mathbb{E}_{\boldsymbol{\mu}_k} \mathbb{X}_{n+1}\mathbf{h}_k \right)\mathbf{h}_k^\top = \sum_{k=1}^{N_\omega} \varpi_k \mathbf{d}_{\boldsymbol{\mu}_k}\mathbf{h}_k^\top$

    \vspace{0.2cm}
    \Statex \hspace{12pt} \textbf{Step 3: convergence check}
    \State \textbf{if} $\max_{1\le k\le N_\omega} \| (\mathbb{X}_{n+1} \!-\! \mathbb{X}_{n}) \mathbf{h}_k \|^2_{\mathbb{M}^\ddagger} / \left( \| \mathbf{c}_{\nu,n} (\boldsymbol{\mu}_k) \|^2_{\mathbb{M}^\dagger} + \| \mathbb{X}_{n+1} \mathbf{h}_k \|^2_{\mathbb{M}^\ddagger}\right) \le \mathrm{tol}^2$
    \State \hspace{14pt} \textbf{break}
\EndFor
\vspace{0.2cm}
\State Set $\mathbb{X} \gets \mathbb{X}_{n+1}$
\end{algorithmic}
\end{algorithm}

\newpage
\subsection{Online Nonlinear System Solution}
\label{sec:online_nonlinear_solve}

Given the lifting matrix $\mathbb{X}$ computed in the offline stage, the online stage evaluates the reduced model for any new parameter value $\boldsymbol{\mu}\in\mathcal{D}$. In particular, we determine the latent coordinates $\mathbf{c}_\nu(\boldsymbol{\mu})\in\mathbb{R}^{N_\nu}$ by solving the nonlinear reduced system~\eqref{eq:nonlinear_rom}, and then reconstruct the high-dimensional coefficient vector via~\eqref{eq:nonlinear_manifold}. Exploiting the affine decomposition~\eqref{eq:affine_decomposition} and the precomputed reduced blocks~\eqref{eq:reduced_blocks}, we assemble, for the queried parameter $\boldsymbol{\mu}$, $\mathbb{A}^\dagger_{\boldsymbol{\mu}}$, $\mathbb{A}^\ddagger_{\boldsymbol{\mu}}$, and $\mathbf{f}^\dagger_{\boldsymbol{\mu}}$. This online assembly is independent of the full-order dimension $N_h$ and scales only with the reduced dimensions and the number of affine terms. 

Specifically, we define the reduced residual
\begin{equation}
\label{eq:online_reduced_residual}
\mathbf{r}_\nu(\mathbf{c}_\nu;\boldsymbol{\mu}) \coloneq \mathbb{A}_{\boldsymbol{\mu}}^\dagger\,\mathbf{c}_\nu + \mathbb{A}_{\boldsymbol{\mu}}^\ddagger \,\mathbb{X}\,\mathbf{h}(\mathbf{c}_\nu) - \mathbf{f}_{\boldsymbol{\mu}}^\dagger \in\mathbb{R}^{N_\nu},
\end{equation}
and seek $\mathbf{c}_\nu(\boldsymbol{\mu})$ such that $\mathbf{r}_\nu(\mathbf{c}_\nu(\boldsymbol{\mu});\boldsymbol{\mu})=\mathbf{0}$. This system can be solved using a Newton (or quasi-Newton) method \cite{Kelley2003, More1980, Powell1970}.  Each nonlinear iteration requires evaluating $\mathbf{h}(\mathbf{c}_\nu)$ (and, if needed, its Jacobian) and solving a dense linear system of size $N_\nu$.

Once $\mathbf{c}_\nu(\boldsymbol{\mu})$ is available, we reconstruct the approximate full-order coefficients and field as
\begin{equation}
\label{eq:online_reconstruction}
\mathbf{c}_\varepsilon(\boldsymbol{\mu}) = \mathbb{V}\,\mathbf{c}_\nu(\boldsymbol{\mu}) + \mathbb{U}\,\mathbb{X}\,\mathbf{h}\!\left(\mathbf{c}_\nu(\boldsymbol{\mu})\right), \qquad u_\varepsilon(\boldsymbol{\mu}) = \sum_{j=1}^{N_h} \big(\mathbf{c}_\varepsilon(\boldsymbol{\mu})\big)_j\,\chi_j .
\end{equation}

If a computable coercivity lower bound $\alpha_a(\boldsymbol{\mu})$ is available (or can be efficiently estimated), then an \emph{a posteriori} error indicator follows from~\eqref{eq:algebraic_error_bound}.

\section{Numerical Experiments}
\label{sec:numerical_experiments}

\subsection{Problem Description: Advection-Diffusion of a Contaminant}
\label{subsec:advection-diffusion}

We consider the steady dispersion of a passive scalar over the square domain
$\Omega \coloneq (-1,1)^2$, modeled by a two-dimensional advection--diffusion equation.
Advection is driven by the incompressible velocity field
\begin{equation}
\boldsymbol{\beta}(x_1,x_2) = 
\begin{bmatrix}
\phantom{-} \sin(\pi x_1)\cos(\pi x_2)\\
-\cos(\pi x_1)\sin(\pi x_2)
\end{bmatrix},
\end{equation}
which generates four Taylor--Green vortices \cite{Taylor1937} over $\Omega$. The parameter $\mu\in\mathcal{D}$, with $\mathcal{D} \coloneq [-4,-2]$, controls the diffusion strength through $e^{\mu}$, which can be interpreted as an inverse Péclet number.

We impose homogeneous Dirichlet conditions on $\Gamma_D \coloneq (-1,1)\times\{-1\}$ and homogeneous Neumann conditions on the remaining part of the boundary, $\Gamma_N \coloneq \partial\Omega\setminus\Gamma_D$. Accordingly, we set $\mathcal{V} \coloneq \{ v\in H^{1}(\Omega) : v|_{\Gamma_D}=0\}$, equipped with the energy inner product and norm
\begin{equation}
(w,v)_{\mathcal{V}} \coloneq \int_{\Omega} \nabla w \cdot \nabla v \, d\Omega, \qquad \|v\|_{\mathcal{V}} \coloneq (v,v)_{\mathcal{V}}^{1/2}.
\end{equation}
For $\mu\in\mathcal{D}$, we define the parametric bilinear form $a(\cdot,\cdot;\mu):\mathcal{V}\times\mathcal{V}\to\mathbb{R}$ as
\begin{equation}
a(w,v;\mu) \coloneq \int_{\Omega} e^{\mu}\,\nabla w\cdot\nabla v \, d\Omega + \int_{\Omega} (\boldsymbol{\beta}\cdot\nabla w)\, v \, d\Omega,
\label{eq:case1_bilinear}
\end{equation}
and seek $u(\mu)\in\mathcal{V}$ such that
\begin{equation}
a(u(\mu),v;\mu) = f(v), \qquad \forall\, v\in\mathcal{V}.
\label{eq:case1_variational}
\end{equation}
The forcing functional is defined by $f(v) \coloneq \int_{\Omega} g\,v\,d\Omega$, with $g(x_1,x_2) = \exp\!\big(- (10x_1 + 1)^2 - (10x_2 - 8)^2\big)$.
%with $g\approx \exp\!\big(- (10x_1 + 1)^2 - (10x_2 - 8)^2)$.

Notice that, with the above choice of $(\cdot,\cdot)_{\mathcal{V}}$, the diffusion term is coercive in the energy norm, while the advection term does not contribute to the coercivity as it is skew-symmetric under the present incompressible flow and boundary conditions. In particular, we obtain the coercivity bound $a(v,v;\mu) \ge e^{\mu}\,\|v\|_{\mathcal{V}}^{2}$ for all $ v\in\mathcal{V}$, so that the coercivity constant is $\alpha_a(\mu)=e^{\mu}$.

\subsubsection{Results with Laplacian Eigenfunction Spaces}
\label{sec:case1_results_spectral}

As a hierarchical $\mathbb{M}$-orthonormal basis, we employ eigenfunctions of the Laplacian on $\Omega$ satisfying the boundary conditions defining $\mathcal{V}$. Specifically, we consider $\{\chi_i\}_{i\ge 1}\subset\mathcal{V}$ such that $-\Delta \chi_i = \lambda_i \chi_i$, and we normalize them so that $(\chi_i,\chi_j)_{\mathcal{V}}=\delta_{ij}$. For the present geometry and boundary partition, these eigenfunctions admit an analytic separated form. In particular, for integers $n_i,m_i\in\mathbb{N}$ with $m_i\le n_i -1$, we take
\begin{equation}
\label{eq:laplacian_eigenfunctions}
\chi_i(x_1,x_2) = \dfrac{ \cos \!\Big( \tfrac{\pi}{2}\,(n_i-m_i-1)\,(x_1+1) \Big)\, \sin \!\Big( \tfrac{\pi}{2}\,(n_i+\tfrac{1}{2})\,(x_2+1) \Big)} {\tfrac{\pi}{2} \sqrt{(n_i-m_i-1)^2+(n_i+\tfrac{1}{2})^2}}.
\end{equation}
The modes are ordered by increasing eigenvalue, i.e., $(n_i-m_i-1)^2+(n_i+\tfrac{1}{2})^2 \le (n_{i+1}-m_{i+1}-1)^2+(n_{i+1}+\tfrac{1}{2})^2$.

\begin{figure}[t]
    \centering
    \includegraphics[width=\linewidth]{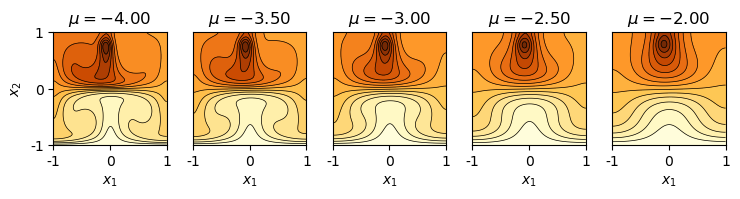}
    \vspace{-0.8cm}
    \caption{Full-order finite element solutions $u_h(\mu)\in\mathcal{V}_h$ for $\mu \in \{-4,-3.5,-3,-2.5,-2\}$.}
    \label{fig:case1_solution}
\end{figure}
\begin{figure}[t]
    \centering
    \includegraphics[width=\linewidth]{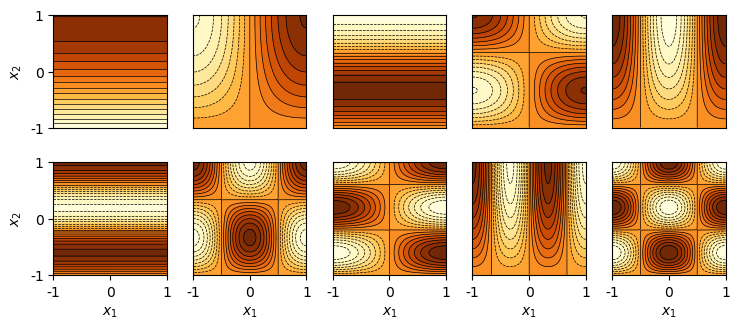}
    \vspace{-0.8cm}
    \caption{First ten $\mathbb{M}$-orthonormal spectral basis functions $\{\chi_i\}_{i=1}^{10}$ (projected on $\mathcal{V}_h$).}
    \vspace{-0.25cm}
    \label{fig:spectral_basis}
\end{figure}
For the numerical experiments, we compute full-order reference solutions in a standard finite element space $\mathcal{V}_h \subset \mathcal{V}$ of dimension $N_h = 40{,}200$, consisting of continuous, piecewise affine functions on a uniform triangulation $\mathcal{T}_h$, i.e., $\mathcal{V}_h=\mathbb{P}^1(\mathcal{T}_h)\cap\mathcal{V}$. These full-order solutions are used \emph{only} as ground truth for visualization and error assessment, and are \emph{not} required by the offline identification procedure described below. All operators and fields (including the analytic spectral modes $\chi_i$) are evaluated on $\mathcal{V}_h$ via projection. Figure~\ref{fig:case1_solution} shows representative solutions $u_h(\mu)$ for five parameter values in $\mathcal{D}$, while Figure~\ref{fig:spectral_basis} reports the first ten spectral basis functions $\{\chi_i\}_{i=1}^{10}$ projected on $\mathcal{V}_h$.

In all experiments we employ a quadratic nonlinear feature map for the latent coordinates \cite{Barnett2022, Geelen2024}. Specifically, for $\mathbf{c}_\nu\in\mathbb{R}^{N_\nu}$ we define
\begin{equation}
\label{eq:quadratic_feature_map}
\mathbf{h}(\mathbf{c}_\nu)
\;\coloneq\;
\mathbb{L}\,\mathrm{vec}\!\left(\mathbf{c}_\nu\,\mathbf{c}_\nu^{\top}\right)
\;\in\;\mathbb{R}^{N_\ell},
\end{equation}
where $\mathrm{vec}(\cdot)$ stacks the entries of a matrix into a vector and $\mathbb{L}\in\mathbb{R}^{N_\ell\times N_\nu^{2}}$ is a fixed selection matrix that extracts only the non-redundant quadratic monomials (e.g., corresponding to the upper-triangular part of $\mathbf{c}_\nu\mathbf{c}_\nu^\top$) and discards duplicates due to symmetry. With this choice, the number of nonlinear features is $N_\ell = N_\nu(N_\nu+1)/2$, and $\mathbf{h}(\mathbf{c}_\nu)$ constitutes the vector of all products $c_{\nu,i}c_{\nu,j}$ with $1\le i\le j\le N_\nu$.

In a first experiment, we set the latent dimension to $N_\nu = 16$ and the lifting dimension to $N_\varepsilon = 500$. For the offline stage, we select a set of $N_\omega = 1{,}001$ parameter samples, uniformly distributed over an \emph{oversampled} interval $\mathcal{D}^+ \coloneq [-4.25,-1.75]$. We emphasize that this set is \emph{not} a training set in the snapshot sense: no full-order solution snapshots are precomputed or stored for these parameters. Instead, the samples are only used to probe the parametric dependence during the snapshot-free training. 

The oversampling $\mathcal{D}^+ \supset \mathcal{D}$ is necessary because we observed a degradation of the identified manifold near the ends of the training interval when sampling only on $\mathcal{D}$; enlarging the interval mitigates these boundary effects and improves accuracy for parameters close to $\mu=-4$ and $\mu=-2$.
A complementary strategy would be to regularize the offline least-squares problem used to identify the lifting matrix $\mathbb{X}$. For example, a Tikhonov-type penalty, such as $\lambda\|\mathbb{X}\|_F^2$, or another problem-dependent regularization term could be added to improve the stability of the learned manifold near the boundary of the parameter domain. This approach has been effective in related time-dependent data-driven manifold approximations, and its extension to the present residual-driven, snapshot-free setting is left for future investigation.

For the test parameters $\mu \in \{ -4.25 + j \times0.25 \}_{j=1}^{9}$, Table~\ref{tab:case1_errors} reports the $\mathcal{V}$-norm errors obtained with: (i) the spectral (linear-subspace) solution in the latent space spanned by the first $N_\nu$ modes; (ii) the nonlinear-manifold solution obtained using the curvature identified by the proposed snapshot-free procedure (training tolerance $\mathrm{tol}=10^{-4}$); and (iii) the best-approximation (projection) error in the lifting space $\mathcal{V}_\varepsilon$ spanned by the first $N_\varepsilon$ modes.

\begin{table}
%\vspace{-0.3cm}
\centering
\caption{$\mathcal{V}$-norm relative errors. Comparison between (i) latent spectral approximation $u_\nu(\mu)$ ($N_\nu=16$), (ii) residual-based nonlinear manifold approximation $u_\varepsilon(\mu)$ ($N_\nu=16$, $N_\varepsilon=500$), and (iii) best approximation in the lifting space $\Pi_{\mathcal{V}_\varepsilon} u_h(\mu)$ ($N_\varepsilon=500$).}
\vspace{0.36cm}
\label{tab:case1_errors}
\begin{tabular}{cccc}
$\mu^\dagger$ & $\dfrac{\|u_h(\mu)-u_{\nu}(\mu)\|_{\mathcal{V}}}{\|  u_h(\mu)\|_{\mathcal{V}}}$ & $\dfrac{\|u_h(\mu)-u_{\varepsilon}(\mu)\|_{\mathcal{V}}}{\|  u_h(\mu)\|_{\mathcal{V}}}$ & $\dfrac{\|u_h(\mu)-\Pi_{\mathcal{V}_\varepsilon} u_h(\mu)\|_{\mathcal{V}}}{\|  u_h(\mu)\|_{\mathcal{V}}}$ \\ 
\midrule
$-4.00$ & $8.3757 \times 10^{-1}$ & $7.0721 \times 10^{-3}$ & $5.5901\times 10^{-3}$ \\
$-3.75$ & $7.9918 \times 10^{-1}$ & $6.1269 \times 10^{-3}$ & $4.2468\times 10^{-3}$ \\
$-3.50$ & $7.5701 \times 10^{-1}$ & $4.9968 \times 10^{-3}$ & $3.6947 \times 10^{-3}$ \\
$-3.25$ & $7.0922 \times 10^{-1}$ & $7.1849 \times 10^{-3}$ & $3.2998 \times 10^{-3}$ \\
$-3.00$ & $6.5449 \times 10^{-1}$ & $6.6490 \times 10^{-3}$ & $2.9652 \times 10^{-3}$ \\
$-2.75$ & $5.9345 \times 10^{-1}$ & $4.6215 \times 10^{-3}$ & $2.6763 \times 10^{-3}$ \\
$-2.50$ & $5.2968 \times 10^{-1}$ & $8.0657 \times 10^{-3}$ & $2.4866 \times 10^{-3}$ \\
$-2.25$ & $4.6876 \times 10^{-1}$ & $8.7753 \times 10^{-3}$ & $2.2142 \times 10^{-3}$ \\
$-2.00$ & $4.1573 \times 10^{-1}$ & $3.9076 \times 10^{-3}$ & $2.0398 \times 10^{-3}$ 
\vspace{-2pt}
\end{tabular}
\end{table}
\begin{table}
\vspace{-0.3cm}
\centering
%\vspace{-1cm}
\caption{$\mathcal{V}$-norm relative residual-based error bounds, normalized by $\alpha_a(\mu)\|u_h(\mu)\|_{\mathcal{V}}$, for $u_\nu(\mu)$, $u_\varepsilon(\mu)$ (residual-based), and $\Pi_{\mathcal{V}_\varepsilon}u_h(\mu)$, similarly to Table~\ref{tab:case1_errors}.}
\vspace{0.36cm}
\label{tab:case1_error_bounds}
\begin{tabular}{cccc}
$\mu^\dagger$ &
$\dfrac{\|r(\,\cdot\,;u_\nu({\mu}), {\mu})\|_{\mathcal{V}'}}{\alpha_a(\mu)\|  u_h(\mu)\|_{\mathcal{V}}}$ & $\dfrac{\|r(\,\cdot\,;u_\varepsilon({\mu}), {\mu})\|_{\mathcal{V}'}}{\alpha_a(\mu)\|  u_h(\mu)\|_{\mathcal{V}}}$ & $\dfrac{\|r(\,\cdot\,;\Pi_{\mathcal{V}_\varepsilon} u_h(\mu), {\mu})\|_{\mathcal{V}'}}{\alpha_a(\mu)\|  u_h(\mu)\|_{\mathcal{V}}}$ \\ 
\midrule
$-4.00$ & $1.9038 \times 10^{-0}$ & $1.0051 \times 10^{-2}$ & $6.6853\times 10^{-3}$ \\
$-3.75$ & $1.5758 \times 10^{-0}$ & $8.0657 \times 10^{-3}$ & $4.7323\times 10^{-3}$ \\
$-3.50$ & $1.3233 \times 10^{-0}$ & $6.3026 \times 10^{-3}$ & $3.9545 \times 10^{-3}$ \\
$-3.25$ & $1.1201 \times 10^{-0}$ & $9.0891 \times 10^{-3}$ & $3.4434 \times 10^{-3}$ \\
$-3.00$ & $9.4633 \times 10^{-1}$ & $7.9746 \times 10^{-3}$ & $3.0445 \times 10^{-3}$ \\
$-2.75$ & $7.9235 \times 10^{-1}$ & $5.3672 \times 10^{-3}$ & $2.7199 \times 10^{-3}$ \\
$-2.50$ & $6.5794 \times 10^{-1}$ & $8.9528 \times 10^{-3}$ & $2.4507 \times 10^{-3}$ \\
$-2.25$ & $5.4691 \times 10^{-1}$ & $9.4988 \times 10^{-3}$ & $2.2276 \times 10^{-3}$ \\
$-2.00$ & $4.6123 \times 10^{-1}$ & $4.1569 \times 10^{-3}$ & $2.0473 \times 10^{-3}$ 
\vspace{10pt}
\end{tabular}
\begin{minipage}{\linewidth}
  \footnotesize
  $^{\dagger}$ The offline training algorithm used the oversampled interval $\mathcal{D}^{+} \coloneq [-4.25,-1.75] $ to mitigate boundary effects near $\mu=-4$ and $\mu=-2$; no full-order snapshots were computed during training. Test parameters are drawn from the original domain
  $\mathcal{D}=[-4,-2]$. \vspace{-0.25cm}
\end{minipage}
\end{table}

Table~\ref{tab:case1_errors} highlights the benefit of enriching the latent spectral truncation with the proposed nonlinear lifting. For all reported parameter values, the manifold approximation $u_\varepsilon(\mu)$ reduces the relative $\mathcal{V}$-error by approximately two orders of magnitude compared with the latent approximation $u_\nu(\mu)$ at the same latent dimension $N_\nu=16$. Moreover, the manifold errors are consistently close to the corresponding best-approximation errors in $\mathcal{V}_\varepsilon$, indicating that the learned nonlinear lifting captures a substantial fraction of the accuracy gain associated with enlarging the approximation space from dimension $N_\nu$ to $N_\varepsilon$, while avoiding full-order solution snapshots during training.

Table~\ref{tab:case1_error_bounds} reports the relative residual-based error bounds under the same testing conditions. The bounds reproduce the same trends observed for the errors: the manifold estimates are approximately two orders of magnitude more accurate than the latent estimates and remain close to the projection benchmarks. In particular, the effectivities are near-unity, confirming that the estimator is sharp.

\begin{figure}[t]
    \centering
    \caption{Spatial error maps for the linear and nonlinear manifold approximations.}
    \label{fig:case1_residual_maps}
    \subfloat[Spatial error maps $u_h(\mu)-u_{\nu}(\mu)$ for the linear spectral approximation with latent dimension $N_\nu=16$, for $\mu\in\{-4.0,-3.5,-3.0,-2.5,-2.0\}$ from left to right.]
    {\includegraphics[width=\linewidth]{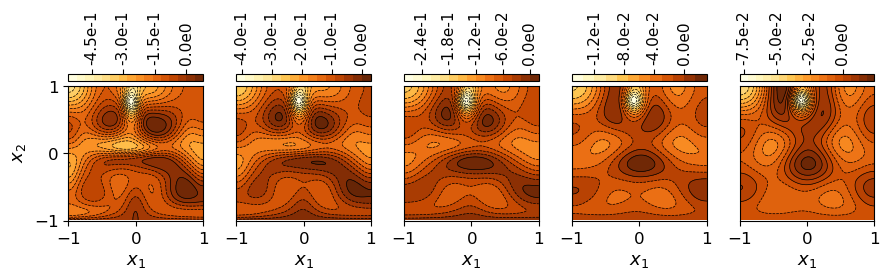}\label{fig:case1_linear_error_maps}}
    \vspace{0.1cm}
    \subfloat[Spatial error maps $u_h(\mu)-u_{\varepsilon}(\mu)$ for the nonlinear manifold approximation with $N_\nu=16$ and $N_\varepsilon=500$, for $\mu\in\{-4.0,-3.5,-3.0,-2.5,-2.0\}$ from left to right.]
    {\includegraphics[width=\linewidth]{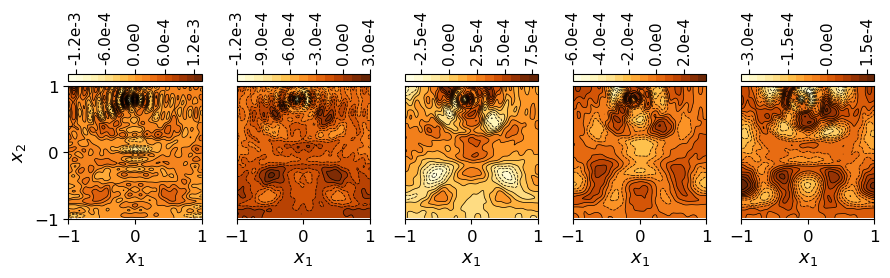}\label{fig:case1_nonlinear_error_maps}}
\end{figure}

To complement the quantitative results in Tables~\ref{tab:case1_errors} and~\ref{tab:case1_error_bounds}, Figures~\ref{fig:case1_residual_maps}a and~\ref{fig:case1_residual_maps}b report respectively the spatial error maps $(u_h(\mu)-u_{\nu} (\mu))/\|u_h(\mu)\|_{L^\infty(\Omega)}$ and $(u_h(\mu)-u_{\varepsilon}(\mu))/\|u_h(\mu)\|_{L^\infty(\Omega)}$ for the five test parameters $\mu\in\{-4.0,\allowbreak-3.5,-3.0,$ $-2.5,-2.0\}$. In both cases, the latent space dimension is fixed to $N_\nu=16$. For the nonlinear manifold approximation, the lifting space dimension is $N_\varepsilon=500$. The linear spectral approximation exhibits large, smooth, coherent error structures across the domain, whereas the nonlinear manifold approximation yields errors that are substantially smaller in magnitude and more spatially localized. This qualitative behavior is consistent with the results in Tables~\ref{tab:case1_errors} and~\ref{tab:case1_error_bounds}, which show that the nonlinear lifting improves the reconstruction accuracy by up to two orders of magnitude compared to linear approximation.

\begin{figure}[t]
    \centering
    \includegraphics[width=0.95\linewidth]{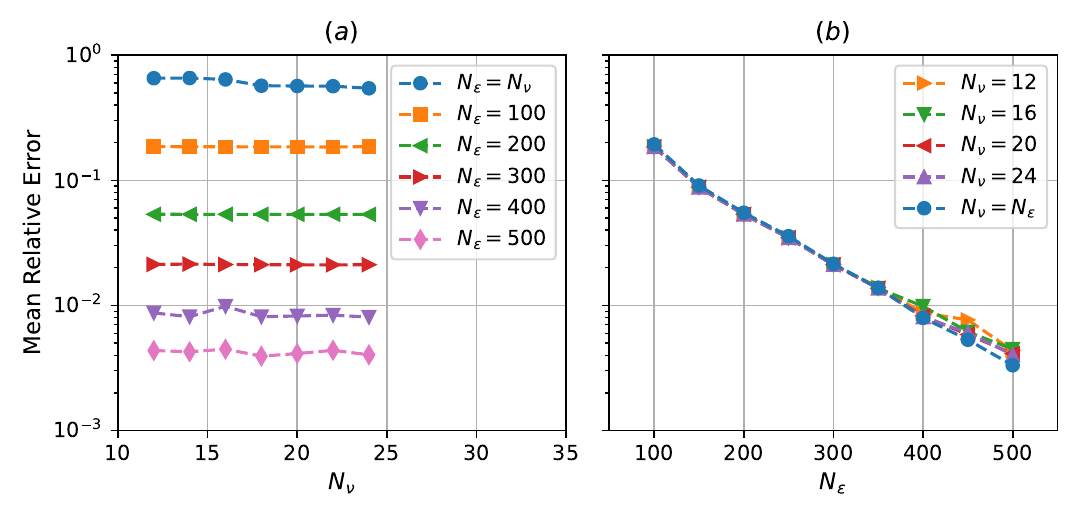}
    \vspace{-0.6cm}
    \caption{Mean relative $\mathcal{V}$-norm error (averaged over the test set) versus latent and lifting dimensions. (a) Error versus $N_\nu$ for fixed $N_\varepsilon \in \{100,200,300,400,500\}$; the linear spectral baseline is $N_\varepsilon=N_\nu$ (blue circles). (b) Error versus $N_\varepsilon$ for fixed $N_\nu \in \{12,16,20,24\}$; the baseline $N_\nu=N_\varepsilon$ is shown in blue with circular markers. All results obtained with $\mathrm{tol}=10^{-4}$.}
    \label{fig:MRE_VS_Nnu_Neps}
\end{figure}
In the next experiment, we study the effect of varying the latent dimension $N_\nu$ and the lifting dimension $N_\varepsilon$. Figure~\ref{fig:MRE_VS_Nnu_Neps} reports the mean relative $\mathcal{V}$-error, averaged over the same test set used in Tables~\ref{tab:case1_errors} and~\ref{tab:case1_error_bounds}, of the nonlinear manifold approximation for $N_\varepsilon\in[100,500]$ and $N_\nu\in[12,24]$, together with the mean relative error of the linear spectral approximation. In both panels, the linear spectral errors are shown in blue with circular markers: in panel~(a) for $N_\nu\in[12,24]$, and in panel~(b) for $N_\nu\in[100,500]$ (i.e., $N_\nu=N_\varepsilon$). Panel~(a) indicates that, over the tested range, the nonlinear manifold error minimally depends on the latent dimension $N_\nu$ and is instead primarily controlled by the lifting dimension $N_\varepsilon$. Panel~(b) confirms that convergence is governed by $N_\varepsilon$: even with the smallest latent space ($N_\nu=12$), increasing $N_\varepsilon$ yields mean errors comparable to those of the linear spectral solver with $N_\nu=N_\varepsilon$.

\begin{figure}
    \centering
    \includegraphics[width=0.6\linewidth]{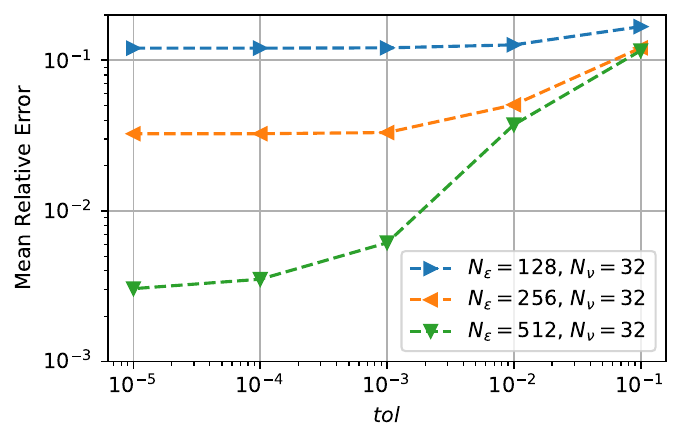}
    \vspace{-0.4cm}
    \caption{Mean relative $\mathcal{V}$-norm error as a function of the training tolerance $\mathrm{tol}$, for nonlinear manifold models trained with fixed latent dimension $N_\nu=32$ and lifting dimensions $N_\varepsilon\in\{128,256,512\}$.}
    \label{fig:tol_vs_error}
\end{figure}

In the third experiment, we investigate the impact of the training tolerance $\mathrm{tol}$. Figure~\ref{fig:tol_vs_error} reports the mean relative $\mathcal{V}$-error over the test set for nonlinear manifold models trained with a fixed latent dimension $N_\nu = 32$ and lifting dimensions $N_\varepsilon \in \{128,256,512\}$, across tolerances $\mathrm{tol} \in \{10^{-5},10^{-4},10^{-3}$, $10^{-2},10^{-1}\}$. The results reveal a clear interaction between $N_\varepsilon$ and $\mathrm{tol}$: increasing $N_\varepsilon$ improves accuracy only when the tolerance is sufficiently small, whereas for larger tolerances training terminates too early to fully exploit the richer lifting space. For the smallest lifting space ($N_\varepsilon=128$), the mean error remains $\mathcal{O}(10^{-1})$ and is largely insensitive to $\mathrm{tol}$. In contrast, for $N_\varepsilon=256$ and $512$, the lowest errors are achieved for $\mathrm{tol} \le 10^{-4}$-$10^{-3}$, while accuracy degrades markedly at larger tolerances. These observations support the choice $\mathrm{tol}=10^{-4}$ adopted in the previous experiments, which is sufficiently small to achieve near-optimal accuracy while controlling training cost.

\subsubsection{Results with Proper Orthogonal Decomposition Spaces}
\label{sec:snapshot_free_pod}

The trends in Figure~\ref{fig:MRE_VS_Nnu_Neps} show that, when using the Laplacian eigenfunctions as hierarchical bases, reducing the approximation error is primarily achieved by increasing the lifting dimension $N_\varepsilon$, while enlarging the latent dimension $N_\nu$ does not yield significant gains. This behavior can be explained by two complementary effects. First, the parameter domain is one-dimensional, hence the solution set $\{u(\mu)\}_{\mu\in\mathcal D}$ is expected to have low intrinsic dimension; beyond a certain threshold, increasing $N_\nu$ does not significantly enrich the manifold's ability to track parametric variability. Second, and more importantly for the present experiment, the achievable accuracy is limited by the approximation properties of the hierarchical basis used for the lifted representation. In our setting, the lifting space $\mathcal V_\varepsilon$ is spanned by Laplacian eigenfunctions, which are not tailored to the advection-dominated structures induced by the Taylor--Green vortices. As a result, the best-approximation error in $\mathcal{V}_\varepsilon$ decays relatively slowly with $N_\varepsilon$, and the nonlinear manifold can only improve the solution to the extent permitted by this fixed representation. This explains why the manifold errors closely track the projection benchmarks in Tables~\ref{tab:case1_errors}--\ref{tab:case1_error_bounds}, and why $N_\varepsilon$ is the dominant accuracy lever in that setting.

\begin{figure}
    \centering
    \includegraphics[width=\linewidth]{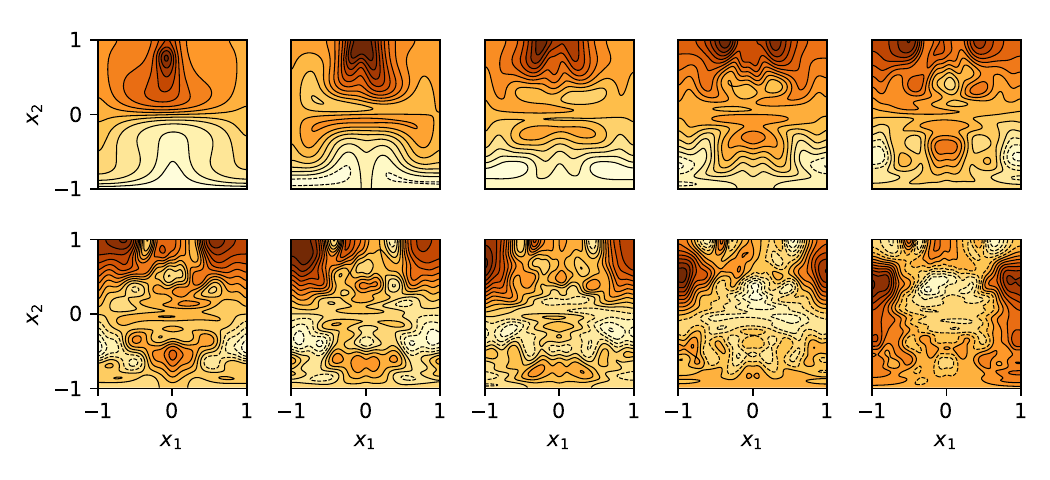}
    \vspace{-0.8cm}
    \caption{First ten $\mathbb{M}$-orthonormal POD basis functions (projected on $\mathcal{V}_h$).}
    \label{fig:POD_basis}
\end{figure}

In practical reduced-order modeling, one would typically replace generic spectral bases with problem-adapted spaces such as proper orthogonal decomposition (POD) bases \cite{Volkwein2008, Sirovich1987}, which often yield significantly faster decay of the projection error. POD spaces are usually constructed from solution snapshots; here, however, we consider a different scenario: we assume that a POD basis is available (e.g., from prior computations or legacy data), but we do not assume access to the snapshot set that generated it. The goal is to assess whether the proposed residual-driven training strategy can still identify an accurate nonlinear lifting \emph{without using solution snapshots during training}. 
This setting should be understood as a controlled test that prescribes a quasi-optimal approximation space in order to isolate the quality of the training procedure from that of the basis, rather than as a fully snapshot-free workflow.
For reference, Figure~\ref{fig:POD_basis} displays the first ten POD basis functions used to build the latent and lifting spaces.

\begin{figure}[b!]
    \centering
    \includegraphics[width=0.55\linewidth]{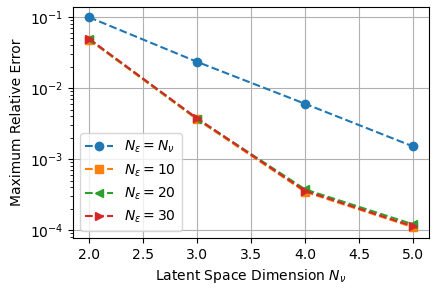}
    \vspace{-0.4cm}
    \caption{Maximum relative $\mathcal V$-norm error over the test set versus latent dimension $N_\nu$ for POD-based spaces. The blue curve corresponds to the linear baseline $N_\varepsilon=N_\nu$, for which the nonlinear lifting term is absent, while the other curves fix $N_\varepsilon\in\{10,20,30\}$.}
    \label{fig:POD_latent}
\end{figure}

\begin{figure}[t!]
    \centering
    \includegraphics[width=0.5\linewidth]{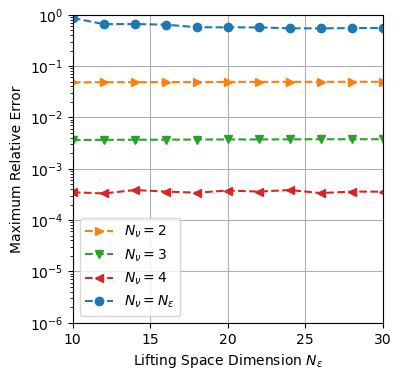}
    \vspace{-0.4cm}
    \caption{Maximum relative $\mathcal V$-norm error over the test set versus lifting dimension $N_\varepsilon$ for POD-based spaces. Curves correspond to fixed $N_\nu\in\{2,3,4\}$, together with the linear baseline $N_\varepsilon=N_\nu$, for which the nonlinear lifting term is absent.}
    \label{fig:POD_lifting}
\end{figure}

Figures~\ref{fig:POD_latent}--\ref{fig:POD_lifting} report the \emph{maximum} relative error over the test set as a function of the latent and lifting dimensions when $\mathbb V$ and $\mathbb U$ are extracted from a POD hierarchy. 
The case $N_\varepsilon=N_\nu$ is included as a linear baseline. Indeed, in this case there are no complementary modes, so the matrix $\mathbb U$ is empty and the nonlinear correction term in \eqref{eq:nonlinear_manifold} vanishes. The approximation therefore reduces to the standard linear truncation on the latent POD space.

Figure~\ref{fig:POD_latent} shows that the error decreases rapidly as the latent dimension $N_\nu$ increases. For the nonlinear approximations with $N_\varepsilon\in\{10,20,30\}$, the curves are nearly indistinguishable. This indicates that, once a moderate lifting space is available, further increasing $N_\varepsilon$ provides little additional benefit. Thus, for POD-based spaces, the accuracy is primarily controlled by the latent dimension $N_\nu$, while the effect of the lifting dimension saturates early.

This conclusion is reinforced by Figure~\ref{fig:POD_lifting}, where the maximum relative error is shown as a function of $N_\varepsilon$ for fixed $N_\nu\in\{2,3,4\}$. Over the tested range $N_\varepsilon\in[10,30]$, the error remains essentially constant for each fixed $N_\nu$, confirming that enlarging the lifting space beyond a moderate size does not lead to a significant improvement in accuracy. The dominant reduction in error is instead obtained by increasing the latent dimension. At the same time, the linear baseline $N_\varepsilon=N_\nu$, for which the nonlinear lifting is absent, remains substantially less accurate. This shows that the nonlinear correction provides a meaningful accuracy gain over the corresponding linear POD truncation, but that this gain is already realized with relatively small lifting dimensions.

Overall, these results indicate that the proposed method can construct nonlinear manifold models whose accuracy is close to that attainable within the prescribed POD-based approximation spaces, once the bases are given, without requiring access to full-order solution snapshots for the identification of the nonlinear lifting. At the same time, the attainable accuracy remains limited by the approximation properties of the selected latent and lifting spaces. In particular, replacing generic spectral bases with POD-based spaces shifts the dominant source of error from the lifting dimension $N_\varepsilon$ to the latent dimension $N_\nu$. As a result, once a moderate lifting space is included, further increasing $N_\varepsilon$ provides only marginal improvement, while increasing $N_\nu$ leads to the main reduction in error.

Importantly, once the prescribed bases and the parameter-independent affine blocks have been precomputed, the repeated lifting updates involve only reduced-dimensional quantities and are independent of the full-order dimension $N_h$. This contrasts with classical snapshot-based nonlinear-manifold workflows, where a dominant cost is often associated with generating, storing, and processing high-fidelity solution snapshots.

\subsection{Problem Description: Compressed Plate with a Circular Cavity}
\label{sec:elasticity}

We consider a second benchmark problem to assess the proposed approach in the vector-valued setting of linear plane-strain elasticity, and to compare the snapshot-free lifting identification of Section~\ref{subsec:optimal_lifting_matrix} against its snapshot-driven counterpart.

The physical system consists of a square plate of side length $L = 48$, containing a circular cavity of radius $R = 8$. Exploiting the geometric and load symmetry, we restrict the computational domain to the first quadrant, $\Omega \coloneq (0, L/2)^{2} \setminus {B(\mathbf{0},R)}$, where $B(\mathbf{0},R) = \{\mathbf{x} \in \mathbb{R}^{2} : |\mathbf{x}| < R\}$ denotes the open disk of radius $R$ centered at the origin. The boundary $\partial\Omega$ is decomposed into five disjoint sections: $\partial\Omega = \Gamma_{L} \cup \Gamma_{B} \cup \Gamma_{T} \cup \Gamma_{H} \cup \Gamma_{R}$, where $\Gamma_{L} \coloneq \{0\} \times (R, L/2)$, $\Gamma_{B} \coloneq (R, L/2) \times \{0\}$ are the left and bottom symmetry boundaries, $\Gamma_{T} \coloneq (0, L/2) \times \{L/2\}$, $\Gamma_R \coloneq \{L/2\}\times(0,L/2)$ are the top loaded face and the free right edge, respectively, and $\Gamma_H \coloneq \partial B(\mathbf 0,R)\cap \partial\Omega$ is the circular hole boundary. The plate is subjected to a uniform vertical compressive displacement $\delta=-1$ on $\Gamma_T$. Symmetry is enforced on $\Gamma_L$ and $\Gamma_B$, while the horizontal displacement is constrained on $\Gamma_T$.
% to model a rigid loading plate.

The parameter $\mu$ denotes the Poisson ratio of the material. Since the volumetric coefficient $\mu/(1-2\mu)$ diverges as $\mu\uparrow1/2$, we exclude the incompressibility limit by a margin $\delta>0$ and consider the closed domain $\mu\in\mathcal{D}\coloneq[-1,\,1/2-\delta]$, taking $\delta=0.02$ in all experiments. The notation $0.5^{-}$ refers to the limiting regime $\mu\uparrow1/2$.

We define the homogeneous function space $\mathcal V_0 \coloneq \{ \mathbf v=(v_1,v_2)\in H^1(\Omega)^2$ with $ v_1=0 \ \text{on } \Gamma_L\cup\Gamma_T,\; v_2=0 \ \text{on } \Gamma_B\cup\Gamma_T \}$. The conditions $v_1=0$ on $\Gamma_L$ and $v_2=0$ on $\Gamma_B$ enforce the geometric symmetry, while the conditions $v_1=v_2=0$ on $\Gamma_T$ correspond to test functions and perturbations satisfying homogeneous essential boundary conditions on the loaded face. We equip $\mathcal V_0$ with the energy inner
product
\begin{equation}
(\mathbf w,\mathbf v)_{\mathcal V} \coloneq \int_\Omega \boldsymbol\varepsilon(\mathbf w):\boldsymbol\varepsilon(\mathbf v) \, d\Omega, \qquad \|\mathbf v\|_{\mathcal V} \coloneq (\mathbf v,\mathbf v)_{\mathcal V}^{1/2},
\label{eq:elast-ip}
\end{equation}
where $\boldsymbol{\varepsilon}(\mathbf{v}) = (\nabla\mathbf{v} + \nabla \mathbf{v}^{\top})/2$ is the symmetric strain tensor. By Korn's inequality \cite{Ciarlet2021} and given the imposed essential constraints, $\|\cdot\|_{\mathcal V}$ defines a norm on $\mathcal V_0$.

For each $\mu\in\mathcal{D}$, the plane-strain bilinear form is defined by
\begin{equation}
a(\mathbf u,\mathbf v;\mu) \coloneq \int_\Omega \boldsymbol\varepsilon(\mathbf u):\boldsymbol\varepsilon(\mathbf v) \, d\Omega + \frac{\mu}{1-2\mu} \int_\Omega \operatorname{tr}\bigl(\boldsymbol\varepsilon(\mathbf u)\bigr) \operatorname{tr}\bigl(\boldsymbol\varepsilon(\mathbf v)\bigr)
\, d\Omega .
\label{eq:elast-bili}
\end{equation}
The first term represents the shear contribution, while the second term is the volumetric contribution. 

The weak problem is posed with test functions in $\mathcal V_0$, after strong imposition of the essential boundary conditions. In the notation of Section~\ref{sec:problem_definition}, this yields a variational problem of the form
\begin{equation}
a\bigl(\mathbf u(\mu),\mathbf v;\mu\bigr) = f(\mathbf v;\mu), \qquad \forall \mathbf v\in\mathcal V_0,
\label{eq:elast-vp}
\end{equation}
where the displacement field $\mathbf u(\mu)$ satisfies the prescribed compressive displacement on $\Gamma_T$ and the symmetry conditions on $\Gamma_L$ and $\Gamma_B$. The functional $f(\cdot;\mu)$ accounts for the elimination of the nonhomogeneous Dirichlet data. No body forces or prescribed
tractions are considered.

We next analyze the coercivity properties of~\eqref{eq:elast-bili}. Since the parameter domain includes negative Poisson ratios, the coefficient $\mu/(1-2\mu)$ is not nonnegative for all $\mu\in\mathcal D$. Nevertheless, the bilinear form is uniformly coercive in the norm~\eqref{eq:elast-ip}. Indeed, for any symmetric second-order tensor in two dimensions, $\left(\operatorname{tr}(\boldsymbol\varepsilon(\mathbf v))\right)^2 \leq 2\,\boldsymbol\varepsilon(\mathbf v):\boldsymbol\varepsilon(\mathbf v)$. Therefore, $\bigl\| \operatorname{tr}(\boldsymbol\varepsilon(\mathbf v)) \bigr\|_{L^2(\Omega)}^{2} \leq 2\|\mathbf v\|_{\mathcal V}^{2}$. It follows that $a(\mathbf v,\mathbf v;\mu) \geq
\alpha_a(\mu)\|\mathbf v\|_{\mathcal V}^{2}$, with
\begin{equation}
\alpha_a(\mu) =
  \begin{cases}
  (1-2\mu)^{-1}, & -1\leq \mu<0,\\
   \;\, 1, & \;\;\,0\leq \mu<1/2.
  \end{cases}
\label{eq:elasticity_coercivity_constant}
\end{equation}
Thus, $\alpha_a(\mu)$ provides a valid coercivity lower bound for the bilinear form on the considered parameter range. In particular, $\alpha_a(\mu)\geq \underline{\alpha}_a=1/3$ for all $\mu\in\mathcal{D}$, so the problem is uniformly coercive. The uniform lower bound $\underline{\alpha}_a$ is used for the residual error bounds and in the snapshot-free training objective. 

Full-order reference solutions are computed using continuous, piecewise affine $\mathbb P^1$ vector fields on a graded triangulation $\mathcal T_h$ of $\Omega$, with local refinement near the hole boundary $\Gamma_H$ to resolve the stress concentration. After enforcing the essential boundary conditions, the homogeneous finite element space is denoted by $\mathcal V_{h,0}\subset\mathcal V_0$, and the resulting number of degrees of freedom is $N_h=30{,}665$.

Since the essential boundary data are parameter-independent, all displacement fields share the same prescribed values on the Dirichlet part of the boundary. We therefore introduce the incompressible-limit displacement
\begin{equation}
\mathbf u_h^\star \coloneq \lim_{\mu\uparrow1/2}\mathbf u_h(\mu),
\label{eq:elast-reference-solution}
\end{equation}
and use it as a parameter-independent reference field. For each $\mu\in\mathcal D$, we decompose
\begin{equation}
\mathbf u_h(\mu) = \mathbf u_h^\star+\mathbf w_h(\mu), \qquad \mathbf w_h(\mu)\in\mathcal V_{h,0}.
\label{eq:elasticity-reference-decomposition}
\end{equation}
The perturbation field $\mathbf w_h(\mu)$ satisfies homogeneous essential boundary conditions and vanishes in the incompressible limit 
\begin{equation*}
\lim_{\mu\uparrow1/2}\mathbf w_h(\mu)=\mathbf 0.
\end{equation*}
The change of variables leads to the following shifted discrete problem: find $\mathbf w_h(\mu) \coloneq \mathbf u_h(\mu)-\mathbf u_h^\star \in \mathcal V_{h,0}$ such that
\begin{equation}
a\bigl(\mathbf w_h(\mu),\mathbf v_h;\mu\bigr) = f(\mathbf v_h;\mu) - a\bigl(\mathbf u_h^\star,\mathbf v_h;\mu\bigr), \quad \forall \mathbf v_h\in\mathcal V_{h,0}.
\label{eq:elast-shifted-vp}
\end{equation}
Thus, the reduced model is not built for the full displacement field $\mathbf u_h(\mu)$, but for the correction $\mathbf w_h(\mu)$ to the incompressible-limit reference field $\mathbf u_h^\star$. For each finite value $\mu<1/2$, the operator in~\eqref{eq:elast-shifted-vp} is the same elastic operator as in the original problem, while the right-hand side is the residual obtained by inserting $\mathbf u_h^\star$ into the equilibrium equation at parameter value $\mu$. 

\begin{figure}[t!]
  \centering
  \includegraphics[width=\textwidth]{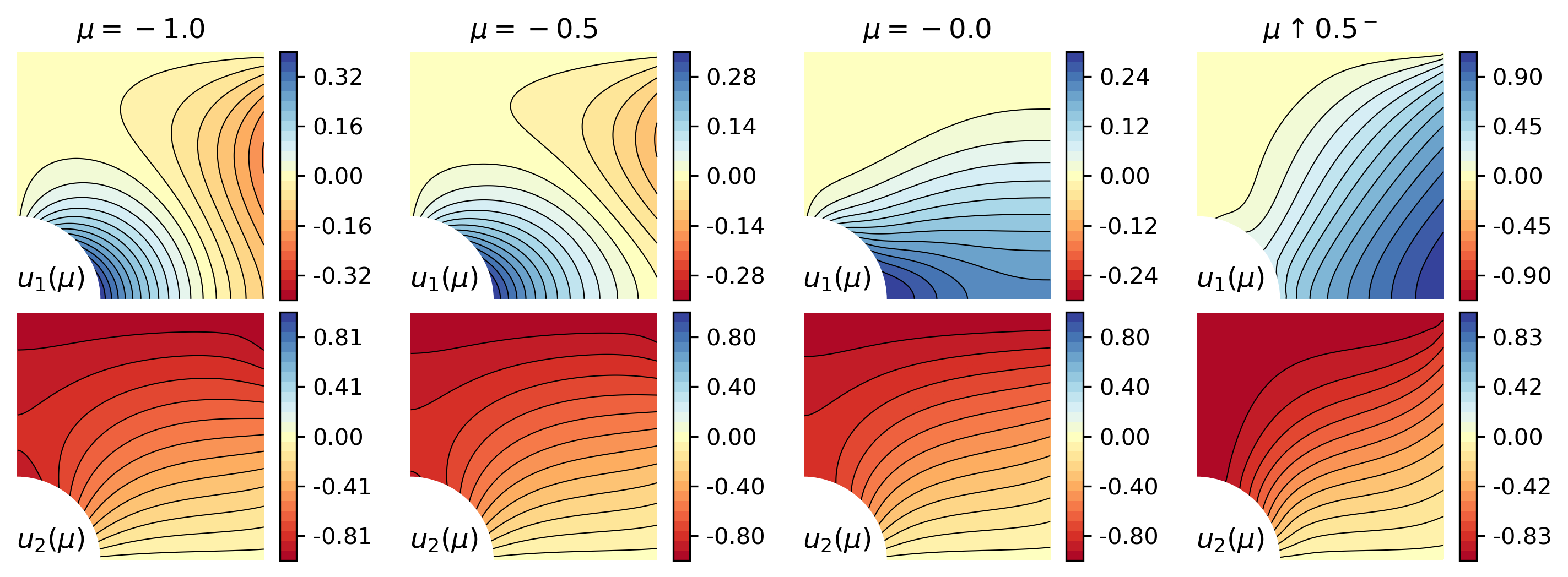}
  \vspace{-0.6cm}
  \caption{Full-order finite element solutions $\mathbf u_h(\mu)$ for the compressed plate problem, for four representative values of $\mu$ ordered
  from $-1.0$ to $0.5^{-}$, left to right. Top row: horizontal component $u_1(\mu)$; bottom row: vertical component $u_2(\mu)$.}
  \label{fig:elast-solutions}
\end{figure}

\begin{figure}[htbp]
  \centering
  \includegraphics[width=\textwidth]{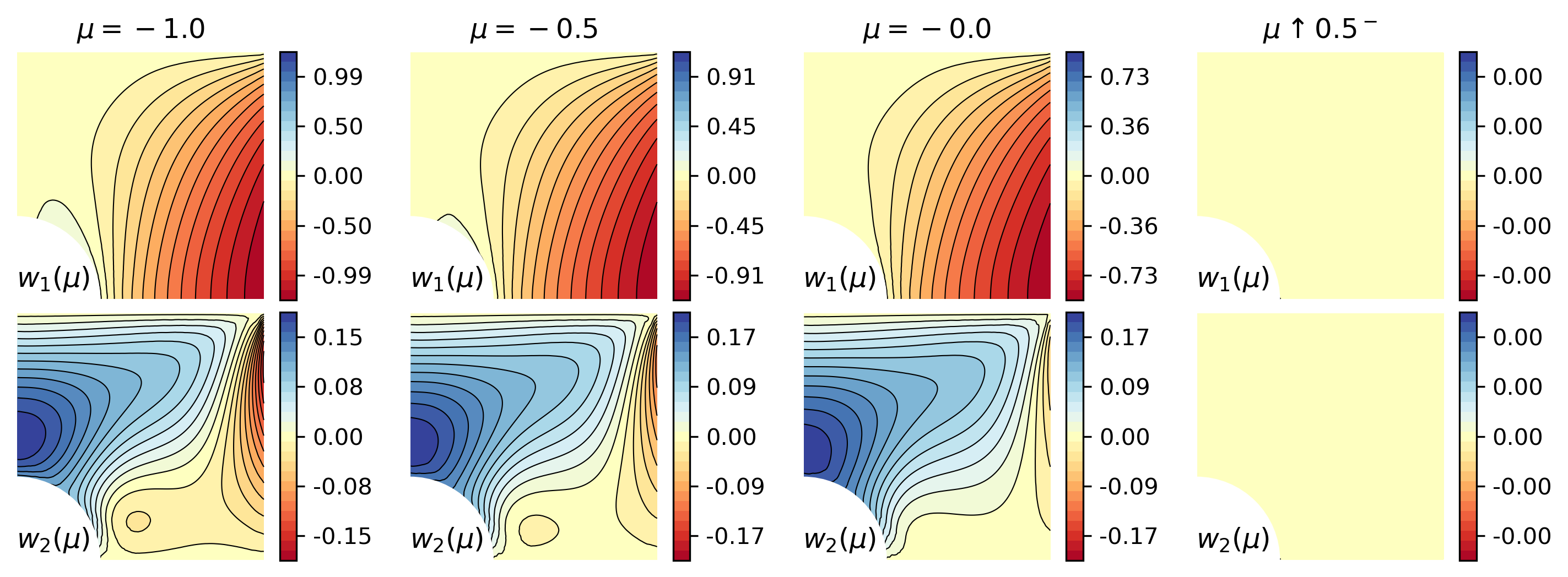}
  \vspace{-0.6cm}
  \caption{Perturbation fields $\mathbf w_h(\mu)=\mathbf u_h(\mu)-\mathbf u_h^\star$ for the compressed plate problem, for the same parameter values as in Figure~\ref{fig:elast-solutions}. Top row: horizontal component $w_1(\mu)$; bottom row: vertical component $w_2(\mu)$. Since $\mathbf w_h(\mu)\to\mathbf 0$ as $\mu\uparrow1/2$, the panel
  corresponding to $0.5^{-}$ is identically zero by construction.}
  \vspace{-0.25cm}
  \label{fig:elast-perturbations}
\end{figure}

Figure~\ref{fig:elast-solutions} reports the $x_1$ and $x_2$-components of the full displacement field $\mathbf u_h(\mu)$ for four representative parameter values ordered from the negative Poisson's ratio regime to the incompressible-limit reference. Figure~\ref{fig:elast-perturbations} shows the corresponding perturbation fields $\mathbf w_h(\mu)$, which constitute the target of the model reduction procedure.

\begin{figure}[htbp]
    \centering
    \includegraphics[width=\linewidth]{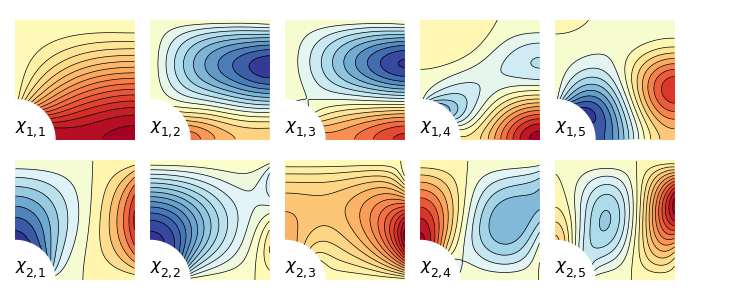}
    \vspace{-0.8cm}
    \caption{First five $\mathbb{M}$-orthonormal spectral basis functions $\{\boldsymbol\chi_i\}_{i=1}^{5}$ for the compressed plate problem. Top row: $x_1$-components; bottom row:
    $x_2$-components.}
    \label{fig:elast-basis}
\end{figure}

As a hierarchical $\mathbb{M}$-orthonormal basis for the latent and lifting spaces, we employ eigenfunctions of the vector-valued strain operator on $\Omega$ with the homogeneous boundary conditions defining $\mathcal V_{h,0}$. Specifically, we consider $\{\boldsymbol\chi_i\}_{i\geq1}\subset\mathcal V_{h,0}$ such that
\begin{equation}
\int_\Omega \boldsymbol\varepsilon(\boldsymbol\chi_i): \boldsymbol\varepsilon(\mathbf v) \, d\Omega = \lambda_i \int_\Omega \boldsymbol\chi_i\cdot\mathbf v \, d\Omega, \quad \forall \mathbf v\in\mathcal V_{h,0}, \quad i=1,2,\ldots,
\label{eq:elast-eig}
\end{equation}
with normalization $(\boldsymbol\chi_i,\boldsymbol\chi_j)_{\mathcal V} = \delta_{ij}$. The modes are ordered by increasing eigenvalue $\lambda_i$. A small value of $\lambda_i$ corresponds to a displacement pattern with small elastic energy relative to its $L^2(\Omega)$ norm; hence the leading modes capture the smoothest admissible displacement patterns. Unlike the scalar Laplacian eigenfunctions of Section~\ref{subsec:advection-diffusion}, these vector-valued modes do not admit a closed-form expression for the present geometry and boundary partition. They are computed numerically using the ARPACK-based sparse symmetric eigensolver \texttt{scipy.sparse.linalg.eigsh}~\cite{2020SciPy}. On our machine, computing the first $260$ basis functions required approximately $12$ seconds. Figure~\ref{fig:elast-basis} displays the $x_1$- and $x_2$-components of the first five spectral basis functions.

\subsubsection{Comparison of Lifting Identification Strategies}
\label{sec:elasticity_data_free_data_based}

We now compare the proposed snapshot-free lifting identification with a snapshot-driven nonlinear manifold construction. In both cases, we use the quadratic nonlinear feature map \eqref{eq:quadratic_feature_map}, with latent dimension $N_\nu=4$ and lifting dimension $N_\varepsilon=256$. As a baseline, we also report the linear spectral approximation $\mathbf u_\nu(\mu)$ obtained in the latent space of dimension $N_\nu$.

The two nonlinear lifting matrices are computed as follows:
\begin{itemize}
  \item the \emph{snapshot-driven} lifting $\mathbb X_{\mathrm{data}}$ is  trained from the full set of high-fidelity snapshots following the data-based nonlinear manifold approach of~\cite{Geelen2024}, with regularization parameter $\lambda_{\mathrm{data}}=10^{-2}$;
  \item the \emph{snapshot-free} lifting $\mathbb X_{\mathrm{free}}$ is computed using Algorithm~\ref{alg:nlmanifold_training_online}, with tolerance $\mathrm{tol}=10^{-3}$, and without computing any full-order solution snapshots during training.
\end{itemize}
All experiments were run on an Apple M2 processor with 8~GB of unified memory. The snapshot-driven workflow used $N_s=1024$ high-fidelity snapshots, computed on a uniform grid over $\mathcal{D}=[-1,\,0.48]$. The snapshot-free training procedure used the same $1024$ parameter values as quadrature samples, thereby isolating the effect of replacing solution snapshots with residual evaluations.

The complete snapshot-driven offline pipeline, comprising snapshot generation and the subsequent regularized least-squares fit of the lifting matrix, required approximately $555$ seconds, with the cost dominated by snapshot generation. In comparison, the proposed snapshot-free training procedure required approximately $29$ seconds, mainly due to the alternating iteration, with a smaller contribution from the assembly of the parameter-independent blocks in~\eqref{eq:offline_blocks}. This figure includes the single full-order solve for the incompressible-limit reference field $\mathbf u_h^\star$, whose cost (less than half a second) is negligible relative to the alternating iteration. The resulting speed-up of about $19\times$ is therefore primarily attributable to avoiding the $1024$ full-order solves required by the snapshot-driven workflow.

These timings cover all auxiliary full-order work across three stages: (i)~constructing $\mathbb{V},\mathbb{U}$ (the strain-eigenfunction eigensolve, ${\approx}12$~s); (ii)~identifying $\mathbb{X}$ (${\approx}29$~s, including the single reference solve); and (iii)~the online reduced solves. None of these uses high-fidelity snapshots. As discussed in Section~\ref{sec:snapshot_free_pod}, only stage~(i) can rely on snapshots, namely when the spaces are POD-based, as in the advection--diffusion experiment.

\begin{figure}[b!]
    \centering
    \includegraphics[width=0.85\linewidth]{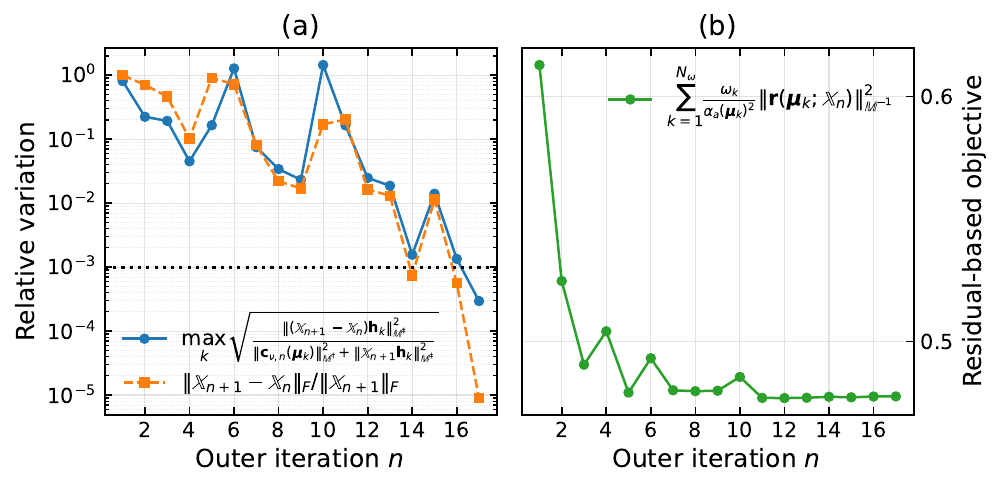}
    \vspace{-0.36cm}
    \caption{Convergence behavior of the alternating training procedure for the compressed-plate problem with $N_\nu=4$, $N_\varepsilon=256$, and $\mathrm{tol}=10^{-3}$. 
    (a) Relative change across outer iterations, measured by the square root of the stopping criterion used in Algorithm~\ref{alg:nlmanifold_training_online}, together with the Frobenius-norm relative change $\|\mathbb X_{n+1}-\mathbb X_n\|_F/\|\mathbb X_{n+1}\|_F$. 
    (b) Evolution of the residual-based training objective.}
    \label{fig:training_history}
\end{figure}
Although a rigorous convergence analysis of the alternating training procedure is beyond the scope of this work, the algorithm is observed to be well behaved in the numerical experiments. Figure~\ref{fig:training_history} reports the training history for the snapshot-free lifting matrix $\mathbb X_{\mathrm{free}}$ in the compressed-plate problem. The residual-based objective is not strictly monotone over the outer iterations, which is consistent with the alternating structure of the algorithm, but it exhibits a clear overall decrease. The convergence indicators in panel~(a) show the same qualitative behavior. In particular, the stopping criterion used in Algorithm~\ref{alg:nlmanifold_training_online} and the simpler Frobenius-norm relative change of the lifting matrix remain close throughout the iteration and cross the prescribed tolerance $\mathrm{tol}=10^{-3}$ within one iteration of each other. Similar behavior was observed for the advection--diffusion benchmark of Section~\ref{subsec:advection-diffusion}, but is not reported here for brevity.

\begin{figure}[b!]
    \centering
    \includegraphics[width=0.9\linewidth]{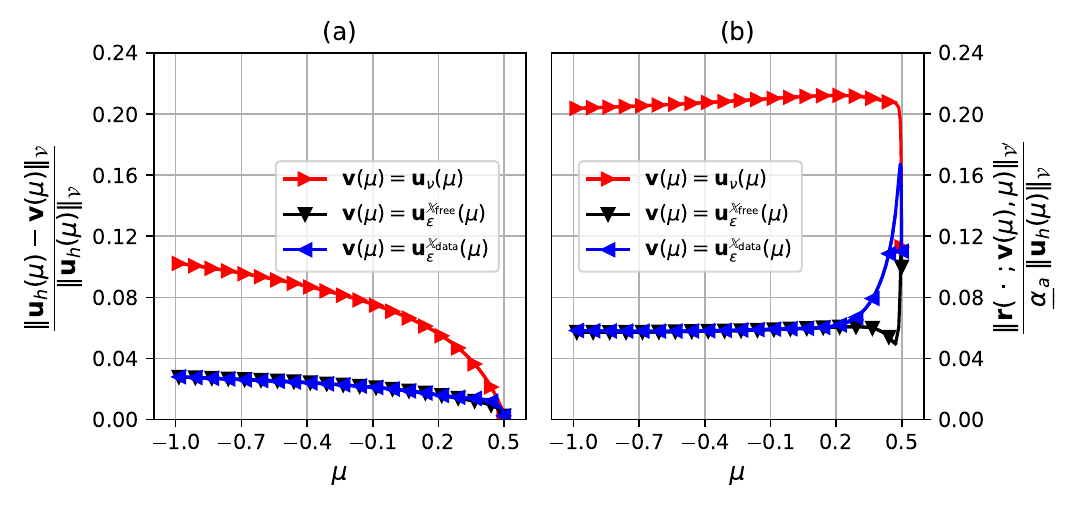}
    \vspace{-0.5cm}
    \caption{Compressed plate problem. Left: relative $\mathcal V$-norm error as a function of the Poisson ratio $\mu$ for the linear spectral approximation, the snapshot-free nonlinear manifold approximation, and the snapshot-driven nonlinear manifold approximation. Right: corresponding residual dual norm normalized by $\underline{\alpha}_a\|\mathbf u_h(\mu)\|_{\mathcal V}$. The latent dimension is $N_\nu=4$ and the lifting dimension is $N_\varepsilon=256$.}
    \label{fig:elast-error-bound-comparison}
\end{figure}
Figure~\ref{fig:elast-error-bound-comparison} reports the results as functions of the Poisson ratio $\mu$. The left panel shows the relative error, while the right panel shows the corresponding residual-based error bound, namely
\begin{equation*}
\frac{\|\mathbf u_h(\mu)-\mathbf v(\mu)\|_{\mathcal{V}}} {\|\mathbf{u}_h(\mu)\|_{\mathcal{V}}} \quad \text{and} \quad \frac{\|r(\,\cdot\,;\mathbf v(\mu),\mu)\|_{\mathcal{V}^\prime}} {\underline{\alpha}_a\|\mathbf{u}_h(\mu)\|_{\mathcal{V}}},
\end{equation*}
respectively, for
$\mathbf v(\mu)\in \{\mathbf u_\nu(\mu), \mathbf u_\varepsilon^{\mathbb X_{\mathrm{free}}}(\mu), \mathbf u_\varepsilon^{\mathbb X_{\mathrm{data}}}(\mu)\}$. 

The nonlinear manifold approximations substantially improve over the linear spectral approximation across most of the parameter range. Near $\mu=-1$, the linear spectral error is approximately $\mathcal O(10^{-1})$, whereas both nonlinear manifold approximations reduce the error to about $\mathcal O(10^{-2})$. As $\mu$ increases toward the incompressible-limit reference, all errors decrease. This behavior is expected since the reduced variable vanishes as $\mu\uparrow1/2$.

The snapshot-driven and snapshot-free nonlinear models yield very similar error cur\-ves. This is the main observation of the experiment: the residual-driven training procedure identifies a lifting that is competitive with the snapshot-driven one, despite not using high-fidelity solution snapshots during training. The agreement is particularly notable because the two liftings are obtained from different objectives: the snapshot-driven model minimizes a reconstruction loss over solution snapshots, whereas the snapshot-free model minimizes a residual-based objective assembled from the affine operator blocks.

The right panel of Figure~\ref{fig:elast-error-bound-comparison} shows that the residual-based bounds reproduce the qualitative separation between the linear and nonlinear approximations over most of the parameter domain. The linear approximation has a significantly larger bound, while the two nonlinear models remain close and substantially lower. Close to the incompressible limit, however, the residual bounds lose sharpness: they increase even though the actual errors decrease toward zero. This behavior is consistent with the growth of the volumetric stiffness as $\mu\uparrow1/2$, which increases the continuity constant of the bilinear form and deteriorates the effectivity of the residual-based $\mathcal V$-norm estimator. Hence, the near-limit rise in the residual bound should be interpreted as a loss of estimator sharpness, rather than as a deterioration of the nonlinear manifold approximation itself. 

We also observe that the snapshot-free model was trained using the parameter-inde\-pendent coercivity lower bound $\underline{\alpha}_a$ rather than the parameter-dependent value $\alpha_a(\mu)$. Since $\underline{\alpha}_a$ is constant, it only rescales the residual objective and therefore drops out of the training algorithm. Despite this conservative choice, the resulting nonlinear manifold exhibits good approximation properties, with a residual-bound that remains nearly uniform over almost the entire parameter domain. By contrast, the residual-based bound associated with the snapshot-driven construction starts to deteriorate slightly earlier as $\mu$ approaches the incompressible limit.

\begin{figure}[htbp]
  \centering
  \includegraphics[width=\textwidth]{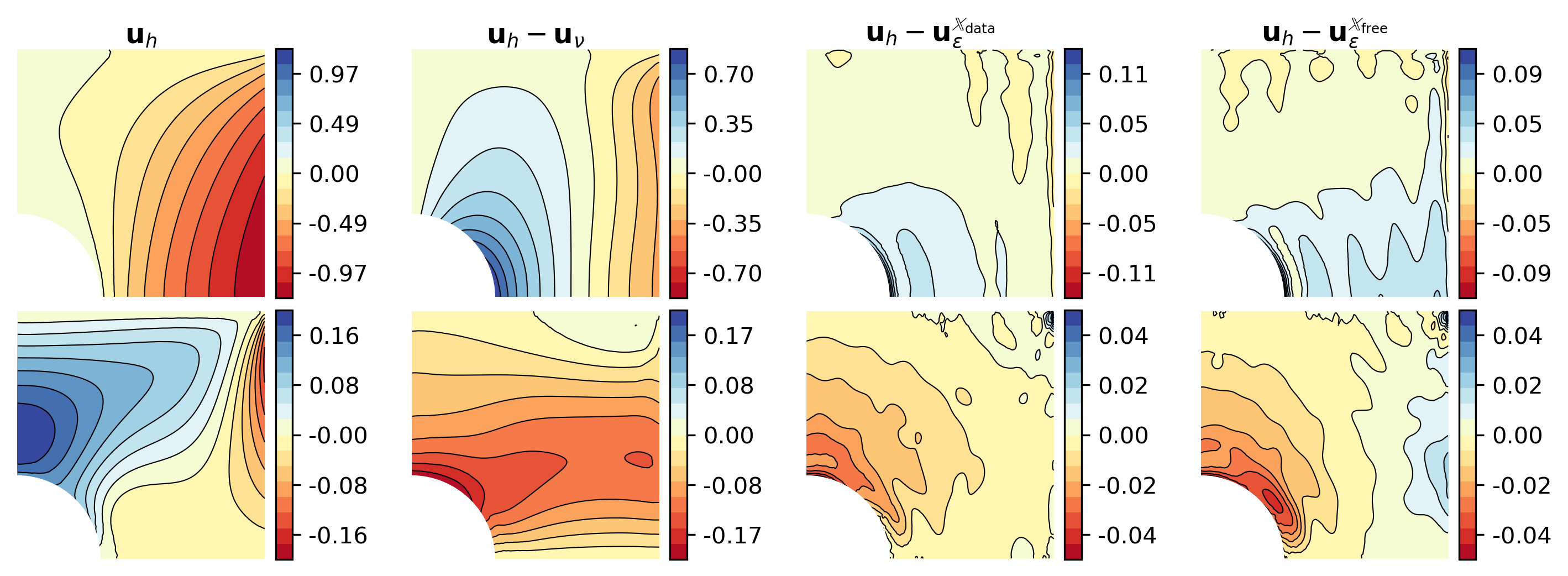}
  \vspace{-0.4cm}
  \caption{Compressed plate problem at $\mu=-1.0$. Top row: $x_1$-component; bottom row: $x_2$-component. From left to right: full-order reference solution $\mathbf u_h(\mu)$; error of the linear spectral approximation $\mathbf u_h(\mu)-\mathbf u_\nu(\mu)$; error of the snapshot-driven nonlinear manifold approximation $\mathbf u_h(\mu)-\mathbf u_\varepsilon^{\mathbb X_{\mathrm{data}}}(\mu)$; error of the snapshot-free nonlinear manifold approximation $\mathbf u_h(\mu)-\mathbf u_\varepsilon^{\mathbb X_{\mathrm{free}}}(\mu)$.}
  \label{fig:elast-errors}
\end{figure}

Finally, Figure~\ref{fig:elast-errors} complements the quantitative results of Figure~\ref{fig:elast-error-bound-comparison} by showing the full-order solution and the corresponding component-wise error fields for the representative parameter value $\mu=-1.0$. The first column reports the full-order displacement $\mathbf u_h(\mu)$, while the remaining columns show the errors associated with the linear spectral approximation, the snapshot-driven nonlinear manifold approximation, and the snapshot-free nonlinear manifold approximation.

The linear spectral approximation exhibits a large, smooth, and spatially coherent error over the whole plate. The error is particularly pronounced near the cavity and along the right portion of the domain. Both nonlinear manifold approximations substantially reduce the error magnitude. In particular, the broad low-frequency error visible in the linear approximation is largely removed. The remaining errors for the nonlinear models are smaller in amplitude and more localized, mainly near the hole boundary, the right edge, and the loaded top boundary. The snapshot-driven and snapshot-free liftings produce comparable error fields, in agreement with the error curves reported in Figure~\ref{fig:elast-error-bound-comparison}.

\section{Conclusions}
\label{sec:conclusions}

We introduced a residual-driven training strategy for non\-lin\-e\-ar-manifold reduced-order models of parametrized coercive linear partial differential equations. The main objective was to identify the nonlinear lifting matrix without relying on high-fidelity solution snapshots. Instead of minimizing an average reconstruction error over precomputed solutions, the proposed method minimizes a computable residual-based upper bound for the state error. The method is therefore snapshot-free \emph{with respect to lifting identification}: no high-fidelity solution snapshots are used to identify the nonlinear lifting matrix. It still requires, however, prescribed latent and lifting spaces, which may themselves be built from snapshots, and access to the affine discrete operators used to assemble the residual objective. It is thus not a fully snapshot-free, end-to-end workflow.

The resulting formulation combines a low-dimensional nonlinear manifold ansatz with standard residual-based \emph{a posteriori} error estimation. For affinely parametrized operators, we derived an offline--online decomposition of the residual objective. After precomputing parameter-independent reduced blocks, the lifting update can be written as a weighted least-squares problem involving only reduced quantities. The resulting alternating algorithm updates the latent coordinates through nonlinear reduced solves and updates the lifting matrix through an action-based solution of the associated normal equations. Once trained, the online stage requires only the solution of a nonlinear system in the latent dimension and the subsequent reconstruction through the learned lifting.

The numerical experiments provide evidence for the accuracy of the proposed snapshot-free training strategy. For the scalar advection--diffusion problem, the re\-si\-du\-al-driven nonlinear manifold substantially improved over the linear spectral approximation, reducing the relative $\mathcal V$-norm error by approximately two orders of magnitude for the tested configuration. The resulting errors were close to the best-approximation errors in the prescribed lifting space, indicating that the residual-driven procedure was able to exploit most of the accuracy available in the enriched approximation space without using solution snapshots during training. Additional tests showed that, when Laplacian eigenfunctions are used as hierarchical bases, the accuracy is primarily controlled by the lifting dimension. In contrast, when POD spaces are prescribed, the error is mainly governed by the latent dimension and the effect of increasing the lifting dimension saturates rapidly.

The compressed-plate elasticity benchmark demonstrated that the approach also extends to vector-valued problems. In this setting, the snapshot-free nonlinear manifold achieved accuracy comparable to a snapshot-driven nonlinear manifold construction, while substantially improving over the linear spectral approximation. The comparison also highlighted the role of the residual-based estimator: over most of the parameter range, the residual bound correctly reflected the separation between the linear and nonlinear approximations, whereas near the incompressible limit the bound became less sharp due to the growth of the continuity constant. This behavior points to the importance of robust residual estimators for nearly incompressible or otherwise parameter-singular regimes. Overall, the results indicate that residual information alone can be sufficient to identify useful nonlinear manifold corrections once appropriate approximation spaces have been selected. This provides a pathway toward nonlinear reduced-order models that retain the accuracy advantages of manifold-based approximations while reducing the dependence on expensive snapshot generation in the offline stage.

Several directions remain open. First, the present work assumes fixed latent and lifting spaces; jointly optimizing these spaces together with the nonlinear lifting could further improve accuracy and reduce the required lifting dimension. Second, a convergence analysis of the alternating training algorithm is still needed, including conditions under which the residual objective admits stable minimizers and the fixed-point iteration converges. Finally, the extension to non-affine, nonlinear, time-dependent, or noncoercive problems would broaden the applicability of the method, but may require empirical interpolation, hyper-reduction, or alternative stability estimates.

%%% ====================================
%%% ====================================
\section*{Acknowledgements}
%%% ====================================
%%% ====================================

The authors used ChatGPT and Claude for language editing and stylistic suggestions. The authors reviewed and are responsible for all scientific content, derivations, code, results, and conclusions.

\bibliographystyle{unsrt}
\bibliography{corrected_references}

\end{document}